\documentclass[12pt]{article}

\usepackage[english]{babel} 
\usepackage[T1]{fontenc} 
\usepackage[utf8]{inputenc} 
\usepackage{latexsym}
\usepackage{euscript}
\usepackage{amsmath,amssymb,MnSymbol,mathrsfs}
\usepackage{amsthm}
\usepackage{enumerate}
\usepackage{enumitem}
\usepackage{longtable}
\usepackage{imakeidx}
\usepackage{graphics}
\usepackage{pictex}
\usepackage{xcolor}
\usepackage[all]{xy}
\usepackage{hyperref}

\usepackage[normalem]{ulem} 

\usepackage{booktabs}

\newcommand {\N}{{\rm I}\hspace{-1mm}{\rm N}}

\newcommand{\fh}{\Rightarrow}
\newcommand{\fD}{{\sf D}}

\newcommand{\set}[1]{\{#1\}}

\newtheorem{definition}{\bf Definition}[section]
\newtheorem{example}[definition]{\bf Example}
\newtheorem{lemma}[definition]{\bf Lemma}
\newtheorem{theorem}[definition]{\bf Theorem}
\newtheorem{remark}[definition]{\bf Remark}
\newtheorem{corollary}[definition]{\bf Corollary}
\renewenvironment{proof}{\noindent\bf Proof \rm}{\hfill$ \square$ \medskip}

\title{Free $T_n$-Heyting algebras }
\author { Juan Manuel Cornejo \and Andr\'es Gallardo \and Paula Helt \and  Luiz Monteiro \and Ignacio Viglizzo }

\makeindex

\begin{document}
\maketitle
\begin{center}
Instituto de Matemática (INMABB), Departamento de Matemática, Universidad Nacional del Sur (UNS) -CONICET, Bahía Blanca, Argentina\\[3mm]
jmcornejo@uns.edu.ar - andres.gallardo@uns.edu.ar - paulasoledadhelt@gmail.com\\ lfmonteiro0510@gmail.com - viglizzo@gmail.com
\end{center}

\thispagestyle{empty}

\begin{abstract}  ${T}_n$-Heyting algebras form a subvariety of Heyting algebras and are characterized by the property that a chain of prime filters in such an algebra has at most $n-1$ elements. We study some basic properties of these algebras and determine the structure and number of elements of the free ${T}_n$-Heyting algebra with a finite number $r$ of free generators, $FT_n(r)$. 
\end{abstract}

\section{Introduction}

Heyting algebras were introduced in 1930 by Arend Heyting to provide an algebraic semantics for intuitionistic logic \cite{heyting30}. These algebras generalize Boolean algebras by weakening the law of excluded middle, thereby capturing the constructive reasoning characteristic of intuitionism. Since their inception, Heyting algebras have played a central role in the study of non-classical logics and have found applications in topology, category theory, and theoretical computer science.

We use Ivo Thomas'  axiom $T_n$ \cite{TI} to characterize the subvariety of Heyting algebras in which chains of prime filters have at most length $n-1$. We call these $T_n$-Heyting algebras or simply $T_n$-algebras (See Definition \ref{DefTn}).

One  particular subclass of \(T_n\)-algebras was studied by Luiz Monteiro in \cite{LM72}.  He studied the free three-valued Heyting algebras (or \(T_3\)-algebras) determining their structure and cardinality when the set of generators is finite. 

All of the $T_n$-algebras are, in particular \textit{linear Heyting algebras}, or $L$-algebras, that is, Heyting algebras satisfying \((x \fh y) \vee(y \fh x)  =  1\). 
Manuel Abad and Luiz Monteiro studied  in \cite{AbMo} the free $L$-algebras. They provided a complete description of the free $L$-algebra with \(n\) generators, including its prime filter structure and a formula for computing its cardinality. An abstract of this work can be found in \cite{Carnielli85meeting}.

Building upon these contributions, the present work determines the structure and cardinality of the free \(T_n\)-Heyting algebra \(FT_n(r)\) with a finite number \(r\) of free generators for arbitrary \(n \geq 2\). Our results here generalize those by L.~Monteiro and agree in the case $n=3$.

\section{Preliminaries} \label{S1}

\begin{definition}\label{D21}
	An algebra $(H,\wedge,\vee,\fh,0,1)$ with arities $(2,2,2,0,0)$ is said to be a Heyting algebra {\rm \cite{AM55,AM59,AM95b}} if the following equations are satisfied:
	\begin{enumerate}[label=\rm{(H\arabic*)}]
		\setcounter{enumi}{-1}
		\item $0 \wedge x  =  0$,\label{propheyting0}
		\item $x \fh x  =  1$,\label{propheyting1}
		\item $(x \fh y) \wedge y  =  y$,\label{propheyting2}
		\item $x \wedge (x \fh y)  =  x \wedge y$,\label{propheyting3}
		\item $x \fh (y \wedge z)  =  (x \fh z) \wedge (x \fh y)$,\label{propheyting4}
		\item $(x \vee y) \fh z  =  (x \fh z) \wedge (y \fh z)$.\label{propheyting5}
	\end{enumerate}
\end{definition}
\medskip

If \(\mathbf{H}\) is a Heyting algebra, then $(H,\wedge,\vee,0,1)$ is a bounded distributive lattice. In Heyting algebras, the intuitionist negation or pseudocomplement is defined by  $\neg x = x \fh 0$. It is well known (see, for example,\cite{AM95b}) that, if \(\mathbf{H}\) is a Heyting algebra, it is also a poset where $x \le y$  if and only if $x = x \wedge y$ and the following properties are verified: 
\begin{enumerate}[label=\rm{(H\arabic*)}]
	\setcounter{enumi}{5}
	\item If $x \le y$, then $y \fh z \le x \fh z$,\label{propheyting6}
	\item $(x \wedge y) \fh z = x \fh (y \fh z)$,\label{propheyting7}
	\item $\neg x \leq x \fh z$,\label{propheyting8}
	\item  $x \fh y = 1$ if and only if $x \le y$,\label{propheyting9}
	\item $\neg x \fh x = \neg \neg x$,\label{propheyting10}
	\item $y \le x \fh y$,\label{propheyting11}
	\item $x \fh (x \fh y) = (x \wedge x) \fh y = x \fh y$,\label{propheyting12}
	\item $1 \fh x = x$.\label{propheyting13}
\end{enumerate}

We recall some further definitions and properties of bounded distributive lattices.

\begin{definition} \label{D12} A subset $F$ of a  bounded distributive lattice \(\mathbf{L}\) is a \emph{filter} if {\rm  (F1)} $1\in F$, {\rm (F2)} if $x \in F, \, x\leq y$  then $y\in F$, and {\rm (F3)} for all $x, y \in F$, $x\land y \in F$.
	We will denote with ${\cal F}(L)$ the family of all filters of \(\mathbf{L}\).
	
	\begin{itemize}
		\item A filter is \emph{proper} if $F\ne L$.
		
		\item If $Y$ is an arbitrary subset of elements of a distributive lattice \(\mathbf{L}\), we denote with $F(Y)$ the \emph{filter generated by} $Y$. This is  the least filter of $L$ containing $Y$.
		
		\item If $F$ is a filter we denote $F(F,a)$ the filter generated by $F\cup \{a\}$.
		
		\item A proper filter $P$ of a distributive lattice \(\mathbf{L}\) is a \emph{prime filter} if $a\lor b\in P$ implies $a\in P$ or $b\in P$. Let $X(L)$ be the set of all prime filters of \(\mathbf{L}\).
		
		\item Maximal, if $F$ is proper, and there is no proper filter that properly contains $F$. We denote with  ${\cal M}(L)$  the family of all maximal filters of \(\mathbf{L}\). 
        
		\item Irreducible, if  $F$ is proper and if  $F=F_1 \cap F_2$ with $F_1,F_2 \in {\cal F}(L)$ then $F=F_1$ or $F=F_2$.  We denote with  ${\cal I}(L)$  the family of all irreducible  filters of \(\mathbf{L}\).
		
		\item Completely irreducible, if  $F$ is proper and if $F=\bigcap\limits_{i\in I} F_i$ with  $F_i \in {\cal F}(L)$  for all $i\in I$ then there exists an index  $i\in I$  such that $F=F_i$. We denote with  ${\cal CI}(L)$  the family of all completely irreducible  filters of \(\mathbf{L}\).
		
		\item Bounded to an  element $f\neq 1$ if it is maximal among those that do not include $f$.
	\end{itemize}
\end{definition}
\medskip

\begin{definition}\label{D14} If \(\mathbf{L}\) is a bounded distributive lattice we denote with ${\bf S}(Q)=\{F\in {\cal F}(L):Q \subseteq F\}$, for all $Q\in{\cal F}(L)$.
It is clear that  ${\bf S}(Q) \neq \emptyset$ and ${\bf S}(Q)$ is a poset when endowed with the relation $\subseteq$.
\end{definition}

\begin{remark}\label{O15} If \(L\) is a distributive lattice, the following are well-known facts: 
	\begin{enumerate}
		\item ${\cal CI}(L) \subseteq {\cal I}(L)$.
		\item ${\cal I}(L) = X(L).$  Therefore, ${\cal CI}(L) \subseteq X(L)$.
		\item If $F \in {\cal F}(L)$ and $F \ne L$, there  exists $M \in {\cal M}(L)$ such that $F \subseteq M$.	
		\item ${\cal M}(L) \subseteq {\cal I}(L)$.
	\end{enumerate}
    Since it follows that all maximal filters are prime, we write $X_M(L)$ for the set of all maximal filters of $L$.
\end{remark}

\begin{theorem} \label{T25a}
	If $H$ is a Heyting algebra, $F \in \mathcal{F}(H)$, $F \neq H$ and  $a \notin F$, then there exists $F_a \in \mathcal{CI}(H)$ such that $F \subseteq F_a$ and $a \notin F_a$.
\end{theorem}

\begin{proof}
	By hypothesis, $F$ is a proper filter.  Let 
	\[
	\mathcal{F}_a = \{ F' \in \mathcal{F}(H)\mid\ F \subseteq F',\ a \notin F'\}.
	\]
	It is easy to check that in $\mathcal{F}_a$ every chain has an upper  bound, so by Zorn's Lemma, there exists a maximal element in $\mathcal{F}_a$, namely $F_a$.
	
	Now we prove that $F_a$ is completely irreducible. If $F_a = \bigcap_{i \in I} F_i$ with $F_i \neq F_a$ for all  $i \in I$, then by the maximality of $F_a$ it follows that  $a \in F_i$ for all $i\in I$, and therefore $a\in F_a$, a contradiction. 
\end{proof}

\begin{lemma}
	\label{L26}
	If $H$ is a Heyting algebra, $F$ a filter and $a\in H$ then the filter generated by $F$ and $a$ is: $$F(F,a) = \{x \in H : a \fh x \in F\} = \{x \in H: \exists\; f\in F, a \wedge f \leq x\}.$$
	If $Y = \{y_1,y_2,\cdots, y_r\} \subseteq H$ then 
	$F(Y) = \left\{ y\in H: \left( \bigwedge \limits_{i=1}^r y_i \right) \fh y = 1 \right\}.$
\end{lemma}

\begin{theorem} \label {T25}{\rm (A. Monteiro, \cite{AM61}, \cite{AM80})}
	If $H$ is a Heyting algebra and $F$ is a filter, then the following conditions are equivalent:
	\begin{enumerate}
		\item $F \in {\cal CI}(H)$, \label{09092026-1}
		\item $F$ is a  filter bound to an element $a \neq 1$, \label{09092026-2}
		\item there exists an element $a\neq 1$ such that  $F$ is a filter  verifying {\rm (i)} $a \notin F$ and \\ {\rm (ii)} $x \fh a \in F$ for every $x \notin F$.\label{09092026-3}
	\end{enumerate} 
\end{theorem}
\begin{proof}
	Let  $F$ be a completely irreducible filter. By definition, $F$ is a  proper filter. For each element  $a \in H \setminus F$, by  Theorem  \ref{T25a}, there exists a filter  $F_a$ containing $F$ and such that $a \notin F_a$. Then:
	\[
	F =\bigcap_{a \in H \setminus F} F_a,
	\]
	and since $F$ is completely irreducible, $F = F_a$ for some $a$.
	
	Assume that $F_a$ is a  filter bounded to an element $a$. If $F_a = \bigcap_{i \in I} F_i$, then for some  $i \in I$, $a \notin F_i$. Also $F_a \subseteq F_i$, but from the definition of bounded filter, we have that  $F_a = F_i$, so $F_a$ is completely irreducible. 

    The equivalence with condition \ref{09092026-3} follows from the characterization of the filter generated by $F$ and $a$ given in Lemma \ref{L26}.
\end{proof}

\begin{lemma}\label{L211} {\rm (A. Monteiro \cite{AM59})}
	If $H$ is a Heyting algebra and $F \in {\cal F}(H)$ then the relation  defined on $H$ by: $$x \equiv_F y \; \; \text{ if  and only if }  \; x \fh y, y \fh x \in F$$ is a congruence of Heyting algebras. 
\end{lemma}

If $x \in H$, we write $C_F(x) = \{y \in H : y \equiv_F x\}$, the equivalence class modulo $F$ containing the element $x$.  We will call the quotient algebra of \(H\) by \(F\), denoted by \(H/F\), the set of all equivalence classes determined by the relation \(\equiv_F\), algebraized in the natural way. The mapping \(C_F \colon H \to H/F\) is called the canonical epimorphism.
from $H$ onto $H/F$.  It is well known that $C_F(1) = F$ and that $C_F(0) = \{x \in H:\neg x \in F\}$. Furthermore: 
\begin{equation}
	\label{ec1}
	C_F(0) \le C_F(x)\le  C_F(1),  \text{ for every } x \in H.
\end{equation}

\section{\texorpdfstring{$T_n$}{Tn}-Heyting algebras}\label{HT}

\subsection{The chains \texorpdfstring{$C_n$}{Cn} and the characterization of  \texorpdfstring{\(T_n\)}{Tn}-algebras}\label{S3}

If $n \in \N$, $n\geq 2$, let $\boldsymbol{C}_n = \left\{a_{j,n} = \dfrac{j}{n-1}: 0 \leq j \leq n-1\right\}$. This set with its natural ordering is a bounded distributive lattice, with first element $0$ and last element $1$. $C_n$ is a chain and a Heyting algebra where the implication is defined by:
\[x \fh y = \begin{cases}
	1 \;{\rm if} \; x \leq y,\cr
	y \;{\rm if} \; x > y.
\end{cases}\]  

\begin{remark}
	\label {O31} 
    Clearly, if  $X \subseteq \mathbf{C}_n$, then $S = \{0,1\} \cup X$ is a Heyting subalgebra of $\mathbf{C}_n$,
	and every Heyting subalgebra of $\mathbf{C}_n$, is isomorphic to $\mathbf{C}_j$ for some  $2 \le j \leq n$.
\end{remark}

\begin{definition}\label{defcadenasdej}
Let $C_n(j) = \{ 0 = a_{0,n}, a_{1,n}, \cdots, a_{j-2,n}, a_{j-1,n} = 1 \}$, where $2 \leq j \le n$, so $C_n(j) \cong C_{j}$, for $2 \leq j \leq n$. $C_n(2) = \{0,1\}$, $C_n(3) = \{0,a_{1,n}=\frac{1}{n-1},1\}, C_n(4) = \{0,a_{1,n},a_{2,n},1\}$ and $C_n(n) = C_n$.	When there is no need to specify the chain \(C_n\), we write \(a_j\) instead of \(a_{j,n}\), for all \(0\le j\le n\).
	
Let $X = \{x_1,x_2,\cdots,x_{j-2}\}\subseteq C_n \setminus \{0,1\}$ where $3 \leq j \leq n$ and  $x_1<x_2<\cdots< x_{j-2}$.  Then $S = \{0=x_0,x_1,x_2,\cdots,x_{j-2},x_{j-1}=1\}$ is a  Heyting subalgebra of $C_n$ with $j$ elements. Furthermore, $S\cong C_n(j)$: it is enough to define $\alpha (x_i)=a_{i,n}$, $0 \le i\leq j-2$ and $\alpha(x_{j-1})=1$. Therefore, if $S$ is a Heyting subalgebra of $C_n$ we can claim that $S\cong C_n(j)$ for some $j$, $2\leq j\leq n$. 
\end{definition}

In his work \cite{TI}, Ivo Thomas advanced the understanding of Dummett's $LC$ logic by developing finite axiomatizations for the Gödel logics associated with finite-valued versions of $LC$. He introduced specific formulas to bridge the gap between the infinite-valued semantics inherent in Dummett's $LC$ and the finite-valued semantics used in Gödel logics. These axioms impose restrictions on the extent of implicational sequences in the following manner: consider well-formed formulas $p_0, p_1, \dots, p_{n-1}$ within the system, and define:
\[\gamma_i = p_i \to p_{i+1} \quad \text{for } i = 0, 1, \dots, n-2,\]
\[\beta_i = \gamma_i \to p_0 = (p_i \to p_{i+1}) \to p_0 \quad \text{for } i = 0, 1, \dots, n-2.\]
The  axiom from Ivo Thomas is then:
\[T_n: \beta_{n-2} \to (\beta_{n-3} \to (\dots \to (\beta_0 \to p_0) \dots )).\]
For $n=2$, this is the Piercean law.

We will now rewrite these axioms as equations defining subvarieties of Heyting algebras.

Let $x_0,x_1,\dots,x_{n-1}$ be variables. We can write $T_1$ for $x_0$, and let for $n\ge 2$:
\[T_n(x_0,\dots,x_{n-1})= ((x_{n-2} \fh x_{n-1}) \fh x_0) \fh T_{n-1}(x_0,\dots,x_{n-2}).\]

\begin{definition} \label{DefTn}
	We say that a Heyting algebra $\bf H$ is a $T_n$-\emph{Heyting algebra} (or $T_n$-\emph{algebra}) if it satisfies the equation:
	\begin{equation}
	    T_n =  1 \tag{$T_n$}
	\end{equation}
\end{definition}

Notice that ${T}_n$-algebras are a particular case of $n$-valued  Hilbert algebras, studied by Luiz Monteiro in \cite{LM77}. Thus, all the results in that work apply to ${T}_n$-algebras.

It follows from the definition that for each $n\ge 2$, $T_n$-algebras form a variety. We will prove that it is the subvariety of Heyting algebras generated by the finite chain $C_n$.  

\begin{remark}\label{O42} Let $p_0,p_1, \dots, p_{n-1}$ be elements of a Heyting algebra $H$. Then:
	\begin{enumerate}
		\item By \ref{propheyting2}, we have that 
		\[T_2(p_0,p_1)\leq T_3(p_0,p_1,p_2)\leq \cdots \leq T_{n-1}(p_0,\ldots,p_{n-2})\leq T_n(p_0,\ldots,p_{n-1}),\]
		 so if $H$ is a  ${T}_n$-algebra, it is also a  ${T}_m$-algebra for all  $m \geq n$.
         
         \item Notice $C_4$ is a ${T}_4$-algebra and not a ${T}_3$-algebra. We calculate in $C_4$:
		\[\left( \left( \dfrac{1}{3} \fh 0 \right) \fh \dfrac{2}{3} \right) \fh \left( \left( \left( \dfrac{2}{3} \fh \frac{1}{3} \right) \fh \dfrac{2}{3} \right) \fh \dfrac{2}{3} \right) =\]
		\[\left( 0 \fh \dfrac{2}{3} \right) \fh \left( \left( \dfrac{1}{3} \fh \dfrac{2}{3}\right) \fh \dfrac{2}{3} \right) =
		1 \fh \left( 1 \fh \dfrac{2}{3} \right) = \dfrac{2}{3} \neq 1.\]
		In a similar way, one can check that for $n\ge 3$, the chain $C_n$ is a $T_n$-algebra and not a $T_{n-1}$-algebra.
		
		\item For \(i=0,\dots,n-2\), we consider $\gamma_i = p_i \fh p_{i+1}$, and $\beta_i = \gamma_i \fh p_0 =$ $ (p_i \fh p_{i+1}) \fh p_0$, we can rewrite the term $T_n$ as: $$T_n=\beta_{n-2}\fh(\beta_{n-3}\fh\cdots\fh (\beta_0\fh p_0)\cdots)){\stackrel{\ref{propheyting7}}{=}} $$ 
		\[\left(\bigwedge\limits_{i=0}^{n-2} \beta_i \right) \fh p_0=\left(\bigwedge\limits_{i=0}^{n-2} (\gamma_i \fh p_0) \right) \fh p_0 {\stackrel{\ref{propheyting5}}{=}} \left(\left( \bigvee\limits_{i=0}^{n-2} \gamma_i\right) \fh p_0 \right) \fh p_0=\]
		\[\left( \left( \bigvee\limits_{i=0}^{n-2} (p_i \fh p_{i+1}) \right) \fh p_0 \right) \fh p_0.\]

        \item From the previous observation, $(T_n)$ can be rewritten as:
\[
  \Bigl(\Bigl(\bigvee_{i=0}^{n-2}(p_i\fh p_{i+1})\Bigr)\fh p_0\Bigr)\fh p_0 = 1.
\]
By (H9) this says $\bigl(\bigvee_{i=0}^{n-2}(p_i\fh p_{i+1})\bigr)\fh p_0\le p_0$,
and since $p_0\le y\fh p_0$ for any $y$ (H11), we obtain the equivalent form:
\begin{equation}
  \Bigl(\bigvee_{i=0}^{n-2}(p_i\fh p_{i+1})\Bigr)\fh p_0 = p_0.
\tag{$T'_n$}  
\end{equation}

		\item  The chains $C_n$ are ${T}_n$-algebras. To see this, let $p_0,p_1,\ldots, p_{n-1}\in C_n$. If for some $i_0$: $p_{i_0} \leq p_{i_0 + 1}$, then $p_{i_0} \fh p_{i_0 +1} = 1$ and we get:
		\[\left( \left( \bigvee\limits_{i=0}^{n-2} (p_i \fh p_{i+1}) \right) \fh p_0 \right) \fh p_0=\left( 1 \fh p_0 \right) \fh p_0 = p_0 \fh p_0 = 1.\]
		
		Assume instead that $p_i > p_{i+1}$ for every $i=0,1,\ldots,n-2$. Thus, $p_0>p_1>\dots>p_{n-1}$ are $n$  different elements in $C_n$, so $p_0=1$ and therefore:
		\[\left( \left( \bigvee\limits_{i=0}^{n-2} (p_i \fh p_{i+1}) \right) \fh 1 \right) \fh 1	= 1.\]

	\end{enumerate}	
\end{remark}

\subsection{Alternative axiomatizations of \texorpdfstring{$T_n$}{Tn}}

Besides Thomas' equation $(T_n)$, seen in the previous section, other axiomatizations of the variety of $T_n$-Heyting algebras have appeared in the literature. Hecht and Katriňák \cite{hecht72equational} proposed in the context of Brouwerian algebras, the equation: 
\[(x_1 \fh x_2 ) \lor (x_2 \fh x_3 ) \lor\ldots\lor (x_n \fh x_{n+1} ) = 1. \tag{$E_n$} \]

More recently, Scirica introduced in \cite{Scirica26Heyting} the equation $(S_n)$ in his characterization of $n$-valued Heyting algebras:
\[\tag{$S_n$} x_1\vee \bigvee\limits_{i=1}^{n-2} (x_i \fh  x_{i+1})\vee \neg x_{n-1} =  1.\]
For $n=2$, we obtain $x \vee \neg x = 1$. Since $x \wedge \neg x  =  0$, this characterizes the variety of Boolean algebras.  For $n=3$, the equation is $(S_3)$:  $x_1 \vee (x_1 \fh x_2) \vee \neg x_2 =  1$ which was also obtained independently by Luiz Monteiro in 2024.

As these equations arise from independent characterizations of the same algebraic phenomenon (the restriction on the length of chains of prime filters) it is expected that they define the same subvariety of Heyting algebras. In this section we study the connections among these axiomatizations. We proved that the Hecht and Katri\v{n}\'ak equation $(E_n)$ implies the Scirica equation $(S_n)$, which in turn implies the Thomas equation $(T_n)$. The proof that $(T_n)$ implies $(E_n)$ was found by Claude Sonnet 4.6; we verified and somewhat simplified it.

\begin{theorem}\label{T35}
For every $n\ge 2$ the equations $(E_n)$, $(S_n)$, and $(T_n)$ are equivalent in Heyting algebras.
\end{theorem}

\begin{proof}
    \begin{description}
      \item[$(E_n)$ implies $(S_n)$] If we replace $x_0$ by $1$ and $x_n$ by $0$ in $(E_n)$, we obtain:
        \begin{align*}
            1 & = (x_0\fh x_1) \vee \bigvee\limits_{i=1}^{{n-2}} (x_i \fh  x_{i+1})  \vee (x_{n-1} \fh x_n) \\
              & = (1\fh x_1) \vee \bigvee\limits_{i=1}^{n-2} (x_i \fh  x_{i+1})  \vee (x_{n-1} \fh 0) \\
              & = x_1 \vee \bigvee\limits_{i=1}^{n-2} (x_i \fh  x_{i+1})  \vee \neg x_{n-1}.
        \end{align*}
\item[$(S_n)$ implies $(T_n)$] We now assume $(S_n)$ holds: 
    \[1 = x_0 \vee \bigvee_{i=0}^{n-3} (x_i \fh x_{i+1}) \vee \neg x_{n-2}\]
    As $x_0 = 1 \fh x_0$, then 
    \[ x_0  = \left( x_0 \vee \bigvee_{i=0}^{n-3} (x_i \fh x_{i+1}) \vee \neg x_{n-2} \right) \fh x_0 \]
    By \ref{propheyting5}, 
    \begin{align*}
        x_0 & = (x_0 \fh x_0) \wedge \left( \left( \bigvee_{i=0}^{n-3} (x_i \fh x_{i+1}) \right) \fh x_0 \right) \wedge (\neg x_{n-2} \fh x_0) \\
            & = \left( \left( \bigvee_{i=0}^{n-3} (x_i \fh x_{i+1}) \right) \fh x_0 \right) \wedge (\neg x_{n-2} \fh x_0).
    \end{align*}
     Again by \ref{propheyting5}, $\displaystyle x_0  = \bigwedge_{i=0}^{n-3} ((x_i \fh x_{i+1}) \fh x_0) \wedge (\neg x_{n-2} \fh x_0)$. Then, since $x_0 \leq ((x_{n-2} \fh x_{n-1}) \fh x_0)$ we get:
     \begin{align*}
     x_0 & = ((x_{n-2} \fh x_{n-1}) \fh x_0) \wedge (\neg x_{n-2} \fh x_0) \wedge \bigwedge_{i=0}^{n-3} ((x_i \fh x_{i+1}) \fh x_0) \\
         & = (((x_{n-2}\fh x_{n-1})\lor\neg x_{n-2})\fh x_0) \wedge \bigwedge_{i=0}^{n-3} ((x_i \fh x_{i+1}) \fh x_0).
    \end{align*}
        By (H8)  \(\neg x_{n-2}\le x_{n-2}\fh x_{n-1}\), so:
    \begin{align*}
        x_0 & = ((x_{n-2}\fh x_{n-1})\fh x_0) \wedge \bigwedge_{i=0}^{n-3} ((x_i \fh x_{i+1}) \fh x_0) \\
            & = \bigwedge_{i=0}^{n-2} ((x_i \fh x_{i+1}) \fh x_0).
    \end{align*}
    Therefore:
    \[1 = x_0 \fh x_0  = \left(\bigwedge_{i=0}^{n-2} ((x_i \fh x_{i+1}) \fh x_0)\right)\fh x_0=\left(\bigwedge_{i=0}^{n-2} \beta_i\right)\fh x_0.\]
    \item[$(T_n)$ implies $(E_n)$]
Let us write $e_i = x_i\fh x_{i+1}$ for $i=1,\dots,n$, and set $e = e_1\vee e_2\vee\cdots\vee e_n$.  We show $e=1$.

Apply $(T_n')$ to the $n$-tuple $(x_2, x_3,\dots,x_{n+1})$:
\[
  \Bigl(\bigvee_{i=0}^{n-2}(x_{i+2}\fh x_{i+3})\Bigr)\fh x_2
  \;=\;
  (e_2\vee e_3\vee\cdots\vee e_n)\fh x_2 \;=\; x_2.
\]
Let $e' = e_2\vee\cdots\vee e_n$.
Since $e'\le e$, (H6) gives $e\fh x_2\le e'\fh x_2 = x_2$.
By (H11), $x_2\le e\fh x_2$.  Hence:
\begin{equation}\label{eq:efhx2}
  e\fh x_2 = x_2.
\end{equation}

Next, apply $(T_n')$ to the $n$-tuple $(e, x_2, x_3,\dots,x_n)$:
\[
  \Bigl((e\fh x_2)\vee(x_2\fh x_3)\vee\cdots\vee(x_{n-1}\fh x_n)\Bigr)\fh e = e,
\]
i.e.\
\[
  \bigl((e\fh x_2)\vee e_2\vee\cdots\vee e_{n-1}\bigr)\fh e = e.
\]
Substituting~\eqref{eq:efhx2}:
\[
  (x_2\vee e_2\vee\cdots\vee e_{n-1})\fh e = e.
\]
Distributing the implication over the join by (H5) we get:
\[
  (x_2\fh e)\wedge(e_2\fh e)\wedge\cdots\wedge(e_{n-1}\fh e) = e.
\]
For each $j\in\{2,\dots,n-1\}$: $e_j\le e$, so $e_j\fh e = 1$ by (H9).
Therefore:
\begin{equation}\label{eq:x2fhe}
  x_2\fh e = e.
\end{equation}
Since $x_2\le x_1\fh x_2=e_1\le e$, we have from \eqref{eq:x2fhe} that $e=x_2\fh e=1$.
    \end{description}
\end{proof}

\subsection{L-algebras}
\begin{definition}\label{DefLineal} {\rm (A. Monteiro, \cite{AM62a}, A. Horn \cite{Hor,Hor1})}  A Heyting algebra $H$, is said to be {\em linear} or an {\em $L$-algebra} if, for every $x,y\in H$, \rm{(HL)} $(x\fh y)\vee (y\fh x)=1$.
\end{definition}

It is well known that the chains $C_n$ are linear algebras. Furthermore, we have:

\begin{lemma}\label{L520}
	If $H$ is a $T_n$-algebra, $n$ a natural number, $n\ge 2$, then $H$ is an L-algebra.
\end{lemma}
\begin{proof}
	To prove that \(H\) is linear, it suffices to show that: if $z\in H$ verifies (i) $x\fh y\le z$ and (ii) $y \fh x\le z$, then $z=1$.
	
	From (i) and \ref{propheyting6} it follows that $1{\stackrel{\ref{propheyting1}}{=}}z\fh z \le (x\fh y)\fh z$, then (iii) $(x\fh y)\fh z=1$. Since $y{\stackrel{\ref{propheyting11}}{\le}}x\fh y{\stackrel{{\rm (i)}}{\le}}z$, then by \ref{propheyting6},
	(iv) $z\fh x\leq y \fh x$. From (iv) and (ii) we have that (v) $z\fh x \leq z$. Hence:
	\[1{\stackrel{\ref{propheyting1}}{=}}z\fh z {\stackrel{\ref{propheyting6}}{\le}} (z\fh x)\fh z.\]
	
	Then (vi) $(z\fh x)\fh z=1$.
	Consider the sequence $p_0=z$, $p_i=x$ for $i$ odd $1\le i \le n-1$, and $p_i=y$ for $i$ even $2\le i \le n-1$, then:
	\begin{itemize}
		\item If $i$ is odd $\beta_{i}=(p_{i}\fh p_{i+1})\fh p_0= (x\fh y)\fh z{\stackrel{{\rm (iii)}}{=}}1$, for $1\le i\le n-2$,
		\item If $i$ is even $\beta_{i}=(p_{i}\fh p_{i+1})\fh p_0= (y\fh x)\fh z{\stackrel{{\rm (ii)}}{=}}1$, for $1\le i\le n-2$, 
		\item $\beta_0=(p_0\fh p_1)\fh p_0= (z\fh x)\fh z{\stackrel{{\rm (vi)}}{=}}1$.
	\end{itemize}
	
	Therefore:
	\[1= \left( \bigwedge\limits_{i=0}^{n-2} \beta_i \right) \fh z=1\fh z=z.\]
\end{proof}

The conclusion of Lemma \ref{L520} also follows from a result of Scirica \cite{Scirica26Heyting}. In the proof of Theorem 3.3.1, Section 3.3, he establishes the same statement by contradiction using the equation $(S_n)$, rather than the equation $(T_n)$ employed here.

\begin{theorem} \label{T53}{\rm  (A.Monteiro, \cite{AM59}, \cite{AM62a}, A. Horn, \cite{Hor, Hor1})} In  Heyting algebras, condition {\rm{(HL)}} is equivalent to:
	\[\tag{H}x\vee y=((x\fh y)\fh y)\wedge ((y\fh x)\fh x).\]
\end{theorem}
\begin{corollary}\label{C54} {\rm  (\cite{AM74a, AM80}, \cite {LM70}.)}
	In linear Heyting algebras:
    \[\neg x \vee x = \neg \neg x \fh x.\]
\end{corollary}
\begin{proof} 
	We have that \((*)\) $x \fh \neg x= x \fh (x \fh 0) {\stackrel{{\rm (H12)}}{=}}x \fh 0=\neg x$. Therefore:
	\begin{align*}
    \neg  x\vee x{\stackrel{\rm(H)}{=}}& ((\neg x\fh x)\fh x)\wedge ((x\fh \neg x)\fh \neg x)\\{\stackrel{{\rm (H10), (*)}}{=}}
	& (\neg \neg x \fh x ) \wedge (\neg x \fh  \neg x)  {\stackrel{{\rm (H1)}}{=}} \neg \neg x \fh x.
    \end{align*}
\end{proof}

\begin{theorem} \label{T56} {\rm  (A. Horn \cite{Hor, Hor1}, A. Monteiro \cite{AM54, AM74, AM74a, AM80}, J. Varlet \cite{Var}.)}
	A Heyting algebra $H$ is an $L$-algebra if and only if one of the following conditions is verified: {\rm (L1)} ${\bf S}(P)$ is a chain, for all  $P\in X(H)$, {\rm (L2)} The family of prime filters which contain a prime filter is a chain.
\end{theorem}

\begin{corollary}\label{C56}
	If $H$ is a $T_n$-algebra, $P \in X(H)$, $F\in {\cal F}(H),$ $F\neq H$, and $P\subset F$  then  $F\in X(H)$.
\end{corollary}

\begin{remark}
If for every $P \in X(H)$ the principal upset ${\bf S}(P)=[P)$ is a chain, then the Hasse diagram of $X(H)$ is a disjoint union of rooted trees; in other words, $X(H)$ is a forest.
\end{remark}

\begin{proof}
	It follows from Lemma \ref{L520} and Theorem  \ref{T56}. 
\end{proof}

\begin{lemma} \label{L56a}{\rm  (A. Monteiro \cite{AM80})}
If $H$ is an  $L$-algebra and $P\in X(H)$ then  $H/P$ is a chain.
\end{lemma}

If $X$ is a finite set, we denote with  $|X|$ its cardinality. We recall one further results about L-algebras:

\begin{lemma} \label{L57} {\rm(\cite{AbMo}, Lemma  3.1 )} If $H$ is an \(L\)-algebra, $G$ a generating set of $H$ with $r$ elements, and $P\in X(H)$ then $|H/P|\leq r+2$. 
\end{lemma}	 

\section{Filters in \texorpdfstring{$T_n$}{Tn}-algebras}

It is well known that filters in Heyting algebras can also be characterized as \textit{deductive systems}: $F$ is a filter of a Heyting algebra if and only if $1\in F$ and:  
\[\tag{D}\label{D} \text{If } x,x\fh y \in F  \text{ then } y \in F. \]

\begin{theorem}\label{T415}
	If $H$ is a  $T_n$-algebra, $n\geq 2$ and $C\in\mathcal{CI}(H)$, 
	then ${\bf S}(C)$ {\rm (see Definition \ref{D14})} has at most $n$ elements. 
\end{theorem}

\begin{proof} Since $C\in\mathcal{CI}(H)$, by Remark \ref{O15}, $C\in X(H)$. By Theorem \ref{T56}, ${\bf S}(C)$ is a chain.  Suppose that there exists a chain of $n+1$ filters containing $C$. We can safely assume that $F_1=C\subset F_2\subset \cdots \subset F_{n}\subset  H=F_{n+1}$. By Corollary \ref{C56} , \(F_i\in X(H)\), for \(1\le i\le n\). Since  $C\in {\cal CI}(H)$ then by Theorem \ref{T25}, $C$ is bounded to an element  $p_0\neq 1$,  $p_0\notin C$. Since $C=F_1\subset F_2$ then $p_0\in F_2$, so $p_0\in F_2\setminus F_1$.

Since each inclusion is strict, for every $i=1,\dots,n-1$ we may choose
\begin{equation}\label{xi}
  p_i\;\in\; F_{i+2}\setminus F_{i+1}.
\end{equation}

Now suppose $\gamma_i=p_i\fh p_{i+1}\in F_{i+2}$.
Since $p_i\in F_{i+2}$, the filter property~(D) gives $p_{i+1}\in F_{i+2}$.
But~\eqref{xi} also gives $p_{i+1}\in F_{i+3}\setminus F_{i+2}$,
so $p_{i+1}\notin F_{i+2}$: contradiction.
Hence $\gamma_i\notin F_{i+2}$ for $i=0,\ldots,n-2$.
Since $C\subseteq F_{i+2}$, every element not in $F_{i+2}$ is not in $C$, so
\begin{equation}\label{gamma-not-C}
  \gamma_i=p_i\fh p_{i+1}\;\notin\; C,\qquad i=0,\dots,n-2.
\end{equation}
 
Let $\beta_i=\gamma_i\fh p_0=(p_i\fh p_{i+1})\fh p_0$. 
By~\eqref{gamma-not-C}, $\gamma_i\notin C=F_1$, so by \ref{T25} we obtain: 
\begin{equation}\label{beta-in-C}
  \beta_i=(p_i\fh p_{i+1})\fh p_0\;\in\; C,\qquad i=0,\dots,n-2.
\end{equation}
 
Since $H$ is a $T_n$-algebra:
\[
  1 = T_n(p_0,p_1,\dots,p_{n-1})
    = \beta_{n-2}\fh\bigl(\beta_{n-3}\fh
      \cdots\fh(\beta_0\fh p_0)\cdots\bigr)
    \;\in\; C,
\]
Applying the filter property~(D) successively, we obtain $p_0\in C$, and this contradicts $p_0\notin C$.
\end{proof}

The previous result also appears in Scirica \cite{Scirica26Heyting}, Theorem 4.3.11. The proof there follows the same strategy by contradiction, but is formulated in terms of the equivalent equation $(S_n)$ rather than $(T_n)$.

\begin{theorem}
	\label{T416} 
	If $H$ is a Heyting algebra such that for all $C \in X(H)$,  ${\bf S}(C)$ is a chain with at most $n$ elements, then $H$ is a $T_n$-algebra ($n\geq 2$).
\end{theorem}	
\begin{proof} 
	Let  $H$ be in the condition of the statement, and let $p_i \in H, 0 \leq i \leq n-1$  be such that:
    \[T_n(p_0, \dots , p_{n-1})\neq 1.\tag{1}\]
	Let $\beta_i =(p_i\fh p_{i+1}) \fh p_0  , 0\leq i \leq n-2$,  $B = \{\beta_0, \beta_{1}, \cdots \beta_{n-2}\}$ and:
   \[\tag{2} F = F(B) = \{x \in H : \left( \bigwedge\limits_{i=0}^{n-2} \beta_i \right) \fh x = 1\}.\]
	Then, from (1) and (2) we have (3) $ p_0 \notin F$; if $p_0 \in F$, by Remark \ref{O42}:
    \[\left( \bigwedge \limits_{i=0}^{n-2} \beta_i \right) \fh p_0 = T_n(p_0, \dots , p_{n-1}) = 1,\] which contradicts (1).
	Therefore, by Theorem \ref{T25a}, $(4)$ there exists $C\in X(H)$ such that $(5)\;  p_0 \notin C$ and $(6) \; F \subseteq C$.\\
	As $B \subseteq F(B) = F$ then $(7) \; \beta_i = (p_i \fh p_{i+1}) \fh p_0 \in F$, for $0 \leq i \leq n-2$; then by $(6)$ and $(7)$ we have that $(8) \; \beta_i \in C$.\\
	From $(8)$ and $(5)$ it follows $(9) \; p_i \fh p_{i+1} \notin C$, for $0 \leq i \leq n-2$: if $p_i \fh p_{i+1} \in C$, as $\beta_i = (p_i \fh p_{i+1}) \fh p_0 \in C$, it would result $p_0 \in C$, which contradicts $(3)$. Also, by \ref{propheyting11}, we have that $(10) \; p_{i+1} \leq p_i \fh p_{i+1}$, so from $(9)$ and $(10)$, we obtain $p_{i+1} \notin C$, for $0 \leq i \leq n-2$ and by $(9)$ we can state that $(11) \; p_{i} \notin C$ for $0 \leq i \leq n-1.$
    
	For $0 \leq i \leq n-1$, let $(12) \; F_i = F(C,p_i) = \{ x \in H : p_i \fh x \in C\}$. Hence, by \ref{propheyting1}, $(13) \; p_i \in F_i$ for $0 \leq i \leq n-1$ and by $(9)$ this gives $(14) \; p_{i+1} \notin F_i$ for $0 \leq i \leq n-1$. Given that $p_{i+1} \in F_i$ consequently, by definition of $F_i$, $p_i \fh p_{i+1} \in F_i$. Which contradicts $(9)$. \\
	Thus, from $(13)$ and $(14)$, $p_{i+1} \in F_{i+1}$ and $p_{i+1} \notin F_i$, for $0 \leq i \leq n-2$, we conclude that $(15) \; F_i \neq F_{i+1}$, $0 \leq i \leq n-2$.\\
	If $F_{i+1} \subseteq F_i$, then as $p_{i+1} \in F_{i+1}$ we obtain $p_{i+1} \in F_{i+1} \subseteq F_i$, this yields $p_i \fh p_{i+1} \in C$ which contradicts $(9)$. Therefore, $F_i \subset F_{i+1}$.
	From $(11)$ and $(12)$ this gives $(16) \; C \subset F_i$, for $0 \leq i \leq n-2$, this shows that $\{C= F_0, F_1, \cdots, F_{n-1},H\} = {\bf S}(C)$ is a chain with $n+1$ elements, absurd.
\end{proof} 

\begin{lemma}\label{L418}
	{\rm (A. Monteiro, \cite{AM80}, Theorem 4.3) }
	In a  Heyting algebra  every proper filter is an intersection of completely irreducible filters.
\end{lemma}

\begin{lemma}\label{L419}
	If $H$ is a \(T_n\)-algebra then ${\cal I}(H) \subseteq {\cal CI}(H).$
\end{lemma}
\begin{proof} Let $I\in {\cal I}(H)$, then $I$ is proper and by Lemma \ref{L418}, $I=\bigcap\limits_{j\in J} I_j$, where \linebreak[4]{$I_j\in {\cal CI}(H)$,} for all $j\in J$.
	Therefore, $I\subseteq I_j$, for all $j\in J$.
	Since $H$ is a \(T_n\)-algebra, then by Lemma \ref{L520}, $H$ is linear. Thus, by Theorem \ref{T56}, $\{I_j\}_{j\in J}$  is a chain. This set cannot have more than $n-1$ elements, otherwise $I_1\subset I_2 \subset \cdots \subset I_n \subset H$
	then the completely irreducible filter $I_1$ would be contained in $n$ filters, which is a contradiction with Theorem \ref{T415}.
	Consequently $\{I_j\}_{j\in J}$ has at most $(n-1)$ elements and since $I=\bigcap\limits_{j=1}^{n-1} I_j$ and $I\in {\cal I}(H)$ then $I=I_j$ for some $j$, $1\leq j \leq n-1$.
	Therefore $I$ is completely irreducible.
\end{proof}

From Remark \ref{O15}, (items 1 and 2) and Lemma \ref{L419} we obtain:

\begin{corollary}\label{C420}                    
If $H$ is a \(T_n\)-algebra, then $X(H)= {\cal I}(H)={\cal CI}(H)$.
\end{corollary}

 By Theorems \ref{T415} and \ref{T416} and Corollary \ref{C420}, we have:
\begin{corollary}\label{C421}
A Heyting algebra $H$ is a $T_n$-algebra if and only if every set of filters containing a $P\in X(H)$ is a chain with at most $n$ elements.
\end{corollary}

\begin{theorem}\label{T424}
	If $H$ is a $T_n$-algebra, and $P\in X(H)$, then:
	\[H/P\cong C_{n}(j)=\{0=a_0,a_1,\cdots,a_{j-2}, 1\}\cong C_{j},\]
	for some $j$, $2\le j\le n$.
\end{theorem}

\begin{proof} Let \( P\in X(H)\). By Theorem \ref{T415} and Corollary \ref{C421}, ${\bf S}(P)$ is a chain, and $2\le |{\bf S}(P)| = j \le n$.
	
	Suppose that:
	\[ P = P_{j-1} \subset P_{j-2} \subset \cdots \subset P_1\subset P_0 = H, \]
	we will show that there are exactly \(j\) equivalence classes (mod. \( P \)):
	\[ E_{j-1} = P_{j-1} ,\quad E_{j-2} = P_{j-2} - P_{j-1},\quad \ldots,\quad E_0 = P_0 - P_{1} = H - P_{1}. \]
	
	Notice that \( P_i,\ 1 \leq i \leq j - 1 \), is a filter bounded to \( x \), for all \( x \in E_{i-1}= P_{i-1} - P_i \), and $1\le i\le j-1$. Otherwise, there is a filter \( J \) such that \( P_i \subset J \) and \( x \notin J \). Then \( J \in {\bf S}(P) \), so $J=P_t$ for $0\le t\le j-1$. Since $P_i\subset J$, we have that $t<i$. On the other hand, since $x\in P_{i-1}$ and $x\notin J$, we cannot have that $t\le i-1$, so $i-1<t<i$, a contradiction. 
	
	\begin{enumerate}[label=\Roman*.]
		\item If \( x, y \in E_{i}, i = 0,1,  \ldots, j-1 \), let us see that \( x \sim y \ (\mathrm{mod}\ P) \), that is, \( x \fh y, y \fh x \in P \).
		
		We know that, for $i=j-1$, $E_{j-1}=P=P_{j-1}$ is an equivalence class mod. $P$.
		
		If $i<j-1$, since \( P_{i+1} \) is bounded to both \( x \) and \( y \), by Theorem \ref{T25}, it follows that \( x \fh y \), \( y \fh x \in P_{i+1} \). 
		\[
		P_{i+1} = P \cup \bigcup_{k=i+1}^{j-2}E_k= P \cup \bigcup_{k=i+1}^{j-2}(P_k- P_{k+1}),
		\]
		
		Assume towards a contradiction that $x\fh y\notin P$. Then, there exists $k_0$, $i+1\le k_0\le j-2$ such that $x\fh y\in E_{k_0}$. As we seen before, $P_{k_0+1}$ is bounded to $x\fh y$, and since $x\notin P_{k_0+1}$, by Theorem \ref{T25}, $x\fh y=x\fh (x\fh y)\in P_{k_0+1}$, a contradiction. Therefore, $x\fh y\in P$ and analogously it can be proved that $y\fh x\in P$.
		
		\item If \( x \in E_{i} \) and \( y \notin E_i \), for \( 0 \leq i \leq j-1 \), let us prove that \( x \not\sim y\ (\mathrm{mod}\ P) \), that is, \( x \fh y \notin P \) or \( y \fh x \notin P \).
		
		We have two cases:

        \begin{enumerate}
            \item If \( x \in P \) and \( y \notin P \) then, if \( x \fh y\in P \), it follows that $y\in P$, a contradiction. Then, $x\fh y\notin P$.
            
            \item If \( x \in E_{i} \) and \( y \notin E_i \), for \( 0 \leq i \leq j-2 \), this yields \( x \in P_i- P_{i+1} \), and $y\in\complement E_i=P_{i+1}\cup\complement P_i$.
            
            We want to show that \( x \fh y \notin P \) or \( y \fh x \notin P \).
            
            Let us assume that $x \fh y \in P$ and $y \fh x \in P$. Then, we have that $x, x \fh y \in P_i$, so $y \in P_{i} \cap \left( P_{i+1} \cup \complement P_{i} \right) = P_{i+1}$. Since $ y \fh x \in P \subset P_{i+1}$, we obtain $ x \in P_{i+1}$, which is a contradiction.
        \end{enumerate}
		
	\end{enumerate}
	
	Conditions I and II show us that \( E_0, E_1, \ldots, E_{j-1} \) are the equivalence classes mod. \(P\). 
	
	We already know that, for all $x\in H$, $P=C_P(1) \ge C_P(x)\ge C_P(0)$. Let us see that $E_{i}\leq E_{i+1}$, for $0\le i\le j-2$. 
	
	Let $x\in E_{i}=P_{i}- P_{i+1}, y\in E_{i+1}=P_{i+1}- P_{i+2}$ (if $i\le j-3$). If $y\in E_{j-1}=P_{j-1}$, as observed previously $C_P(x)\le C_P(1)=P$.
	
	Suppose that $x\fh y\notin P$. From $x,y\notin P$ and Lemma \ref{L520} we deduce that $H$ is an \(L\)-algebra, so $x\fh y\in P$ or $y\fh x\in P$. Therefore, $y\fh x\in P\subseteq P_{i+1}$ and since $y\in P_{i+1}$, it follows that $x\in P_{i+1}$, a contradiction. 	Hence, $x\fh y\in P$, so $E_{i}\fh E_{i+1}= P$.
\end{proof}

\begin{definition} \label{defkernel}
    Given a homomorphism $h: H\to H'$ of Heyting algebras, its \emph{kernel} is:\[\mathrm{Ker}(h)=\{x\in H : h(x)=1\}.\]
\end{definition}

\begin{lemma}\label{lemakernelfiltro}
    If \(H\) is a Heyting algebra, \(\mathrm{Ker}(h)\in\mathcal{F}(H)\). If \(h\) is a Heyting homomorphism, then by Lemma~\ref{L211}, \(\mathrm{Ker}(h)=\{1\}\) if and only if \(h\) is injective.
\end{lemma}

Recall that, for a class $\mathcal{K}$ of algebras, $\mathbb{H}(\mathcal{K})$,
$\mathbb{S}(\mathcal{K})$ and $\mathbb{P}(\mathcal{K})$ denote, respectively,
the classes of all homomorphic images, subalgebras, and direct products of
members of $\mathcal{K}$, while $\mathbb{V}(\mathcal{K})$ denotes the
\emph{variety generated by $\mathcal{K}$}, that is, the smallest variety
containing $\mathcal{K}$; by Birkhoff's theorem this coincides with
$\mathbb{HSP}(\mathcal{K})$.

\begin{theorem}\label{ThmVCn}
For every $n\geq 2$, the variety \(\mathcal{T}_n\) of $T_n$-algebras coincides with the
variety generated by the chain $\mathbf{C}_n$, that is:
\[
\mathbb{V}(\mathcal{T}_n) \;=\; \mathbb{V}(\mathbf{C}_n)\;=\;\mathbb{HSP}(\mathbf{C}_n).
\]
\end{theorem}

\begin{proof}
By Remark~\ref{O42}, $\mathbf{C}_n$ is itself a $T_n$-algebra, so
$\mathbb{V}(\mathbf{C}_n)\subseteq\mathbb{V}(\mathcal{T}_n)$.

For the converse inclusion, let $H$ be a $T_n$-algebra. Since $\{1\}$ is a
proper filter of $H$, Lemma~\ref{L418} guarantees that $\displaystyle\{1\} = \bigcap_{j\in J} I_j$,
where each $I_j$ is a completely irreducible filter of $H$. By Corollary~\ref{C420}, since $H$ is a $T_n$-algebra we have $X(H)={\cal I}(H)={\cal CI}(H)$, so every $I_j$ is in fact a \emph{prime} filter of $H$.
Hence $\displaystyle \{1\} = \bigcap_{P\in X(H)} P$. 

For each $P\in X(H)$, Lemma~\ref{L211} shows that the relation $\equiv_{P}$ associated to $P$ is a congruence, giving rise to the quotient algebra
$H/P$ and the canonical epimorphism $C_{P}\colon H\to H/P$. Consider
the map:
\[
\varphi: H \longrightarrow \prod_{P\in X(H)} H/P, \qquad
\varphi(x) = \big(C_{P}(x)\big)_{P\in X(H)}.
\]
This is a homomorphism of Heyting algebras, and it is injective: it suffices to show that $\mathrm{Ker}(\varphi)=\{1\}$: if $\varphi(x)=\varphi(1)$ then $x\equiv_P 1$ for every $P\in X(H)$, so $\displaystyle x\in\bigcap_{P\in X(H)} P=\{1\}$, whence $x=1$. By Lemma \ref{lemakernelfiltro}, $\mathrm{Ker}(\varphi)=\{1\}$ implies $\varphi$ is injective, so $H$ is isomorphic to its image under $\varphi$.

By Theorem~\ref{T424}, for every $P\in X(H)$ there exists $2\le j\le n$ such that $H/P \cong \mathbf{C}_n(j)$, that is, each quotient $H/P$ is (isomorphic to) a Heyting subalgebra of $\mathbf{C}_n$. Consequently $H$ embeds into a product of subalgebras of $\mathbf{C}_n$, which shows that $H \in \mathbb{SP}(\mathbf{C}_n) \subseteq \mathbb{HSP}(\mathbf{C}_n)$.

Since $H$ was an arbitrary $T_n$-algebra, $\mathbb{V}(\mathcal{T}_n)\subseteq \mathbb{HSP}(\mathbf{C}_n)$. By Birkhoff's theorem \cite{Birkhoff1935structure}, $\mathbb{HSP}(\mathbf{C}_n)=\mathbb{V}(\mathbf{C}_n)$,
the variety generated by $\mathbf{C}_n$. Combining both inclusions, we get $\mathbb{V}(\mathcal{T}_n) =\mathbb{V}(\mathbf{C}_n)$.
\end{proof}

\section{Free Algebras} \label{S5}

If $\boldsymbol{A}$ is a Heyting algebra, and $X\subseteq A$, we denote with $HS_A(X)$ the Heyting subalgebra  of $\boldsymbol{A}$ generated by $X$. If the algebra $\boldsymbol{A}$ is a chain and  $X\subseteq A$ then $HS_A(X) = X\cup \{0,1\}$. When it is clear from the context, we simply write \(HS(X)\),  omitting the algebra \(\boldsymbol{A}\).

The notion of free algebra is defined in the usual way, that is:
\begin{definition}
	Given a cardinal $\beta>0$ we say that a $T_n$-Heyting algebra $F$ has a set of $\beta$ free generators, if:
	\begin{enumerate}
		\item $F$ contains a subset $G$ of cardinality $\beta$ such that $HS_F(G) =F$.
		\item Every application $f:G \to H$, where $H$ is any $T_n$-Heyting algebra, can be extended to a unique  
        Heyting algebra homomorphism $\bar{f}$, from $F$ to $H$.
	\end{enumerate}
	
	Under these conditions we say that $G$ is a set of free generators of $F$ and a $T_n$-Heyting algebra  is said to be \emph{free} if it has a set of free generators. To show explicitly the cardinal $\beta$ we write $F=F({\beta})$.
\end{definition}

Since the notion of $T_n$-Heyting algebra is defined through identities, a result of G. Birkhoff, \cite{Birkhoff1935structure},
ensures that there exists the free $T_n$-Heyting algebra $F({\beta})$  with a free set of generators $G=\{g_i\}_{i\in I}$ of cardinality $\beta$, and it is unique up to isomorphisms. If $r$ is a natural number, $r\geq 1$, we are going to determine the $T_n$-Heyting algebra  with a set of $r$ free generators $G=\{g_1,g_2,\ldots,g_r\}$, and we will denote it with $FT_n(r)$.

\medskip
It is well known that:
\begin{lemma}
	\label{L61}
	If ${H}$ and ${H'}$ are Heyting algebras, $h: H \to H'$ an homomorphism, and $G\subseteq H$ such that $HS_H(G) = H$, then $HS_{H'}(h(G)) =h(H)$. That is, if $G$ generates $H$ then $h(G)$ generates $h(H)$. If $h$ is an epimorphism, then $HS_{H'}(h(G))=H'$.
\end{lemma}

\begin{lemma}\label{LkerXH}
Let $H$ be a Heyting algebra and let $h: H \to \mathbf{C}_i$, $2\le i\le n$,
be an epimorphism. Then $\mathrm{Ker}(h)=h^{-1}(\{1\})\in X(H)$.
\end{lemma}

\begin{lemma} \label{Lhom1}
	If $H,H_1,H_2$ are Heyting algebras, $h_i$ an epimorphism from $H$ to $H_i$, $i=1,2$ and $Ker (h_1) \subseteq Ker (h_2)$, then there exists an epimorphism $h$ from $H_1$ to $H_2$ such that $h\circ h_1=h_2$. If $Ker (h_1)= Ker ( h_2)$ then $H_1$ is isomorphic to $H_2$.
	\end{lemma}

\begin{remark}\label{Rdef}
	From now on, for the purpose of improving the notation, we will write $X(n,r)$ instead of $X(FT_n(r))$, and $X_M(n,r)$ denotes $X_M(FT_n(r))$. Also, when there is no need to specify the algebra \(\mathbf{A}\), we will write \(HS(G)\) instead of \(HS_A(G)\).
\end{remark}

A. Horn {\rm \cite{Hor1}} proved that $FL(r)$, the free L-algebra with a set of $r$ generators, is finite, and he determines its number of elements by logical methods. M.~Abad and L.~Monteiro {\rm \cite{AbMo}, page 9}, proved that  $FL(r)$ is finite by an algebraic method.

\medskip
 
\begin{lemma} \label{L52} $FT_n(r)$ is a  finite $T_n$-algebra.
\end{lemma}
\begin{proof} $|FL(r)|\geq |FT_n(r)|$. Let $G'=\{g'_1,g'_2,\cdots, g'_r\}$ be a set of free generators of   $FL(r)$ {\rm (\cite{AbMo}, \cite{Hor1})}, and  $G=\{g_1,g_2,\cdots, g_r\}$ a set of free generators of  $FT_n(r)$. Hence, in  particular, {\rm  (1)} $HS(G')= FL(r)$ and {\rm  (2)} $HS(G)= FT_n(r)$. Let  $f:G'\to G$ defined by $f(g'_i)= g_i$, $1\leq i \leq r$, since $FT_n(r)$ is a $T_n$-algebra then by Lemma  \ref{L520}, $FT_n(r)$  is an L-algebra. Consequently, $f$ can be extended to a  homomorphism $\bar{f}$ from $FL(r)$ to $FT_n(r)$, in particular, (3) $\bar{f}(G')= G$. By Lemma \ref{L61}, we have that (4) $\bar{f}(HS(G'))=HS(\bar{f}(G'))=HS(f(G'))$. Therefore:
$$\bar{f}(FL(r)){\stackrel{{\rm (1)}}{=}} \bar{f}(HS(G'))  {\stackrel{{\rm (4)}}{=}}HS(\bar{f}(G')){\stackrel{{\rm (3)}}{=}}$$
$$HS(G){\stackrel{{\rm (2)}}{=}}FT_n(r).$$
So $\bar{f}$ is an epimorphism, this implies that $|FL(r)|\geq |FT_n(r)|$.
\end{proof}	
\medskip

Then, since \(FT_{n}(r)\) is finite, it is sufficient to find the set \(X(n,r)\) of its prime filters to determine it. 

\begin{definition} \label{defF_jnr}
Let $F_j(n,r)=\{f\in C_n^G: HS_{C_n}(f(G))= C_n(j)\}$, $2\leq j \leq n$.  It is clear that every $F_j(n,r)$ is nonempty and these sets are pairwise disjoint. Let $\displaystyle F(n,r)=\bigcup_{j=2}^n F_j(n,r)$ then if $f\in F(n,r)$, there exists a unique $j$ such that $f\in  F_j(n,r)$. 
\end{definition}

If $P\in X(n,r)$, then by Theorem \ref{T424}:

\[FT_{n}(r)/P\cong C_{n}(j)=\{0=a_0,a_1, \cdots,a_{j-2}, 1\},\text{ for some } 2 \leq j\le n.\]

\begin{definition}
	Let $X_j(n,r)=\{ P\in X(n,r):FT_{n}(r)/P\cong C_n(j)\},  2 \leq j\leq n$. 
	
	It is clear that $X(n,r)= \bigcup\limits_{j=2}^{n} X_j(n,r)$ and $X_j(n,r)\cap X_k(n,r)=\emptyset$ for $j\neq k$ , $2\leq j,k\leq n$.
\end{definition}

\begin{remark}\label{remarkidentificacioncocientesencadenas}
	If $P\in X_j(n,r)$, we can see the homomorphism $h_P: FT_n(r)\to FT_n(r)/P$ as one having codomain $C_n(j)$. Following the notation in Definition \ref{defcadenasdej} and Theorem \ref{T424}, we can write    \(h_P(x)=a_{t,n}\) if and only if, $x\in E_t$, for $0\le t\le j-2$ and $h_P(x)=1$ if and only if, $x\in P=E_{j-1}$. From now on, we will use this identification, by writing $FT_n(r)/P=C_n(j)$. 
\end{remark}

\begin{lemma}\label{L611} There is a bijection $\psi_j$ between the sets $F_j(n,r)$ and $X_j(n,r)$ ($2\leq j\leq n$).
\end{lemma}
\begin{proof} Since $j\le n$ we have that $F_j(n,r)\neq \emptyset$. Each $f\in F_j(n,r)$, can be extended to a unique homomorphism $\bar{f}$ from $FT_{n}(r)$ to $C_n(j)$. By Lemma \ref{L61}, (*) $\bar{f}(HS(G))= \bar{f}(FT_{n}(r))$. Since
$$C_n(j)= HS_{C_n}(f(G))=HS_{C_n}(\bar{f}(G))=\bar{f}(HS(G)){\stackrel{{\rm (*)}}{=}}\bar{f}(FT_{n}(r)),$$ 
it follows that $\bar{f}$ is an epimorphism from $FT_n(r)$ over $C_n(j)$. By Theorem \ref{T424}, $P_f=Ker(\bar{f})\in X_j(n,r)$. 

We define for $j$, $2\leq j \leq n$ a function $\psi_j: F_j(n,r) \to X_j(n,r)$ by $\psi_j(f)= Ker(\bar{f})$, where $f\in F_j(n,r)$.
	
Let us prove that $\psi_j$ is surjective. Let $P \in X_j(n,r)$ and $h$ the natural homomorphism from $FT_n(r)$ over $FT_n(r)/P= C_n(j)$ and $f=h|_{G}$ the restriction of $h$ to the set $G$, we conclude $HS(f(G))= HS(h(G))=h(HS(G))= h(FT_n(r))= C_n(j)$. Therefore $f\in F_j(n,r)$. Let $\overline{f}$ the extension of $f$. Since $\overline{f}|_G=f=h|_{G}$, this yields $\overline{f}=h$ and therefore $\psi_j(f) = Ker(\overline{f})=Ker(h)=P$.
	
Finally, we now show that $\psi_j$ is injective. If $\psi_j(f)=\psi_j(g)$, that is, $Ker(\overline{f})=Ker(\overline{g})$ and then by Lemma \ref{Lhom1}, there exists an isomorphism $h:C_n(j)\to C_n(j)$ such that $h\circ\overline{f}=\overline{g}$. Since the only automorphism over a finite chain is the identity function, we have that $\overline{f}=\overline{g}$ and therefore $f=g$.
\end{proof}

\begin{remark} \label{O612} Since $F_j(n,r)$, $2\leq j \leq n$ are finite sets, then by Lemma \ref{L611}, the sets $X_j(n,r)$ for $2\leq j \leq n$ are finite, and hence $X(n,r)$ is a finite set. This is another way of showing that $FT_n(r)$ is a finite algebra. 
\end{remark}

\begin{corollary}\label{biyeccion}
	There exists a bijection $\psi:F(n,r)\to X(n,r)$.
\end{corollary}
\begin{proof}
	It is enough to define $\psi(f)=\psi_j(f)$ for $f\in F_j(n,r)$. We will also write \(P_f=\psi(f)\).
\end{proof}

\begin{definition}  If $Z$ is a poset we say that $b\in Z$ covers $a\in Z$ if $ a < b$, but $a  <z < b$ for no $z\in Z.$ 
\end{definition}  

\begin{lemma}\label{covers}
If $P,Q \in X(n,r)$, $P$ covers $Q$ if and only if the following conditions are verified:
	\begin{enumerate}
		\item $Q \subset P$, \label{covers 1}
		\item there exists $j$, $2\leq j \leq n-1$ such that $P\in X_j (n,r)$, and  $Q \in X_{j+1}(n,r)$. \label{covers 2}
	\end{enumerate}
\end{lemma}

\begin{proof}
\textit{Necessary condition}: If $P$ covers $Q$, then clearly $Q \subset P$ and (1) $Q\subset F \subset P$ for no $F\in X(n,r)$. Notice also that, if $Q  \subset F \subset P$ was true for any filter $F$, since $Q$ is prime, $F$ would be prime as well.

Since $P\in X(n,r)$, then $FT_n(r)/P = C_n(j)$ for some $j$ with  $2\leq j \leq n$. If $j=n$, ${\bf S}(Q)$  would have $n+1$ elements, which is a contradiction with Theorem \ref{T424}. Therefore, $FT_n(r)/P = C_n(j)$ with $2\leq j\leq n-1$, so from (1) it follows that $|{\bf S}(Q)|=j+1$ so   $FT_n(r)/Q = C_n(j+1)$ and $Q\in X_{j+1}(n,r)$.

\textit{Sufficient condition}: Let $P,Q \in X(n,r)$ such that conditions \ref{covers 1} and \ref{covers 2} are verified. In particular,  $Q\subset P$. If we suppose there is $F \in X(n,r)$ such that $Q\subset F\subset P$, then if $|{\bf S}(P)|=j$ we have that $|{\bf S}(Q)|\ge j+2$, which contradicts condition 2.
\end{proof}

For each function $f\in F(n,r)$, that is, $f:G\to C_n$, we think of the $r$-tuple $(f_1,\ldots,f_r)$, where $f_i=f(g_i)$ for $1\le i \le r$. Then, since by Corollary \ref{biyeccion} each prime filter in $FT_n(r)$ is uniquely assigned to one of these functions, we can denote the filter as $P_{(f_1,\ldots,f_r)}$. Next, we want to determine the order in $X(n,r)$ based on the $r$-tuples of $C_n^G$.

\begin{theorem} \label{orderFilters}
Let $f,h \in F(n,r)$. Then, for $\psi(f)=P_f$ to cover $\psi(h)=P_h$ it is necessary and sufficient that the following conditions are verified:
	\begin{itemize}
		\item[{\rm (i)}] $f \in F_j(n,r)$, $h \in F_{j+1}(n,r)$,  for some $j$,  $2\leq j \leq n-1$, 
		\item[{\rm (ii)}] $f(g)=a_t$ if and only if $h(g)= a_t$ , $0\leq t \leq j-2$,
		\item[{\rm (iii)}] $f(g)=1$ if and only if $h(g)= 1$ or $h(g)= a_{j-1},$	
		\item[{\rm (iv)}] There exists $g\in G$ such that $f(g)\equal 1$ and $h(g)\equal a_{j-1}$.
	\end{itemize}
\end{theorem}
\begin{proof} 
	\textit{Necessary condition.}  
	Assume that  $P_f$ covers $P_h$. Then by Lemma \ref{covers}, we have that $P_h\subset P_f$ and there exists $j, 2\le j\le n-1$ such that $P_f\in X_j(n,r)$ and $P_h\in X_{j+1}(n,r)$. Using the inverse of $\psi_j$ from Lemma \ref{L611}, this implies that $f\in F_j(n,r)$ and $h\in F_{j+1}(n,r)$, so (i) is verified.
	
	Since $HS(f(G)) = C_n(j)$, then $f(G) \cup \{0,1\} = C_n(j)$. If $j=2$ then $f(G) \subseteq \{0,1\}$, otherwise $a_1, a_2, \dots, a_{j-2} \in f(G)$. 
	Similarly, $a_1, a_2, \dots, a_{j-1} \in h(G)$.
	Since
	\[P_{j} = P_h \subset P_f = P_{j-1} \subset P_{j-2} \subset \dots \subset P_1 \subset P_0 = FT_n(r),\] we have:
	\[\bar{f}(x) = \begin{cases}
		1 & \text{if } x \in P_f \\
		a_t & \text{if } x \in P_t - P_{t+1}, \quad 0 \leq t \leq j-2,
	\end{cases}\]
	and
	\[\bar{h}(x) = \begin{cases}
		1 & \text{if } x \in P_h \\
		a_t & \text{if } x \in P_t - P_{t+1}, \quad 0 \leq t \leq j-1.
	\end{cases}\]
	
	Then, 
	\begin{itemize}
		\item \( \bar{f}(x) = \bar{h}(x) = a_t \) for $0 \leq t \leq j-2$ if and only if $x \in P_t - P_{t+1}$,
		\item \( \bar{f}(x) = \bar{h}(x) = 1 \) if and only if $x \in P_h$,
		\item \( \bar{f}(x) = 1 \) and \( \bar{h}(x) = a_{j-1} \) if and only if $x \in P_f - P_h$.
	\end{itemize}
	
	In particular, if $g \in G$:
	\begin{itemize}
		\item[(ii)] $f(g) = a_t$ if and only if $h(g) = a_t$, for $0 \leq t \leq j-2$.
		\item[(iii)] $f(g) = 1$ if and only if $h(g) = 1$ or $h(g) = a_{j-1}$.
	\end{itemize}
	
	In addition:
	\begin{itemize}
		\item[(iv)] There exists $g \in G$ such that $f(g) \neq h(g)$. Indeed, if $f(g) = h(g)$ for all $g \in G$, then $f = h$ and therefore $P_f = P_h$.
	\end{itemize}
	
	\textit{Sufficient condition.} Let $f, h \in F(n,r)$ satisfying (i), (ii), (iii) and (iv).
	
	From (i), $f \in F_j(n,r)$ and $h \in F_{j+1}(n,r)$  so we have:
	\[FT_n(r)/P_f = C_{n}(j), \quad FT_n(r)/P_h = C_{n}(j+1), \quad 2 \leq j \leq n-1\]
	
	By Lemma \ref{covers},  we only need to prove that $P_h \subset P_f$.
	
	Consider:
	\[
	P_f= P_{j-1} \subset P_{j-2} \subset \dots \subset P_1 \subset P_0 = FT_n(r)\]
	and
	\[P_h = Q_{j} \subset Q_{j-1} \subset \dots \subset Q_1 \subset Q_0 = FT_n(r)\]
	the chains of prime filters containing $P_f$ and $P_h$ respectively, and consider the following sets:
	
	\[D_j = Q_{j} \cap P_{j-1}=P_h\cap P_f,\]
	\[ D_{j-1} = (Q_{j-1} - Q_{j}) \cap P_{j-1}\]
	\[ D_{t} = (Q_{t} - Q_{t+1}) \cap (P_t-P_{t+1}), \text{ for } 0\le t\le j-2.\]
	Then we have:
	\begin{itemize}
		\item $z \in D_{j}$ if and only if $\overline{f}(z) = \overline{h}(z)=1$ .
		\item $z \in D_{j-1}$ if and only if $\overline{h}(z) = a_{j-1}$ and $\overline{f}(z) = 1$.
		\item $z \in D_t$ if and only if $\overline{f}(z) = a_t=\overline{h}(z)$, $0 \leq t \le j-2$.
	\end{itemize} 
	
Observe that $D_0 = (Q_0-Q_1) \cap (P_0-P_1) = \complement(Q_1 \cup P_1)$. 
We note that $D_0$ is an ideal because the complement of any prime filter is a prime ideal---defined as a subset containing $0$ that is downward closed and closed under the join operation. Consequently, $D_0$ is the intersection of two ideals, which is itself an ideal, ensuring that $0 \in D_0 \subseteq S$.
	
We have that $D_{j}$ is a filter, because it is an intersection of filters.

The sets $D_t$, $0 \le t \leq j-2$ are nonempty: since $a_t\in \overline{f}(FT_n(r))=HS(f(G))=f(G)\cup\set{0,1}$, then there exists $g\in G$ such that $f(g)=a_t$. By condition (ii), we also have $h(g)=a_t$ so $g\in D_t$.

The set $D_{j-1}$ is also nonempty: Indeed, from (iv), there exists $g \in G$ such that $f(g) \neq h(g)$, and then we have, from (ii) and (iii), that $\overline{h}(g) = a_{j-1}$ and $\overline{f}(g) = 1$, that is, $g \in D_{j-1}$.

It is clear that the sets $D_t$, $0 \le t \le j$, are pairwise disjoint. Observe that $D_{j} \cup D_{j-1} = Q_{j-1} \cap P_{j-1}$ and so it is a filter.

Let us prove that the set $S = \bigcup\limits_{i=0}^{j} D_i$ is a subalgebra of $FT_n(r)$. It is clear that $1\in D_j\subseteq S$, which is a filter and $0\in D_0\subseteq S$ which is an ideal, then (i) $0,1 \in S$.

In  the following tables, $0\leq t\le q\leq j-2$:
\begin{center}
    \begin{longtable}{||l|c|c|c|c|c|c|c||}
        \hline\hline
                 $x \in$                & $y\in$ & \rule{0pt}{3ex} $\overline{h}(x\vee y)$ & $\overline{f}(x\vee y)$ & $x \vee y \in $ & $\overline{h}(x\wedge y)$ &$\overline{f}(x\wedge y)$ & $x \wedge y \in $ \\ \hline\hline
        $D_j$                                        &            $D_j$             &                      1                       &                      1                       &                $D_j$                 &                       1                        &                       1                        &                  $D_j$                  \\ \hline
        $D_{j}$                                      &          $D_{j-1}$           &                      1                       &                      1                       &                $D_j$                 &                   $a_{j-1}$                    &                       1                        &                $D_{j-1}$                \\ \hline
        $D_{j}$                                      &           $D_{t}$            &                      1                       &                      1                       &                $D_j$                 &                     $a_t$                      &                     $a_t$                      &                  $D_t$                  \\ \hline
        $D_{j-1}$                                    &          $D_{j-1}$           &                  $a_{j-1}$                   &                      1                       &              $D_{j-1}$               &                   $a_{j-1}$                    &                       1                        &                $D_{j-1}$                \\ \hline
        $D_{j-1}$                                    &           $D_{t}$            &                  $a_{j-1}$                   &                      1                       &              $D_{j-1}$               &                     $a_t$                      &                     $a_t$                      &                  $D_t$                  \\ \hline
        $D_{t}$                                      &           $D_{t}$            &                    $a_t$                     &                    $a_t$                     &                $D_t$                 &                     $a_t$                      &                     $a_t$                      &                  $D_t$                  \\ \hline
        $D_{t}$                                      &           $D_{q}$            &                    $a_q$                     &                    $a_q$                     &                $D_q$                 &                     $a_t$                      &                     $a_t$                      &                  $D_t$                  \\ \hline\hline
\end{longtable} \pagebreak[4]
    \begin{longtable}{||l|c|c|c|c|c|c|c||}
        \hline\hline
        $x \in$                & $y\in$ & \rule{0pt}{3ex} $\overline{h}(x\fh y)$  & $\overline{f}(x\fh y)$  & $x \fh y \in $ &  $\overline{h}(y\fh x)$   & $\overline{f}(y\fh x)$   &  $y \fh x \in $   \\
        \hline\hline
        $D_j$                                        &            $D_j$             &                      1                       &                      1                       &                $D_j$                 &                       1                        &                       1                        &                  $D_j$                  \\ \hline
        $D_{j}$                                      &          $D_{j-1}$           &                  $a_{j-1}$                   &                      1                       &              $D_{j-1}$               &                       1                        &                       1                        &                  $D_j$                  \\ \hline
        $D_{j}$                                      &           $D_{t}$            &                    $a_t$                     &                    $a_t$                     &                $D_t$                 &                       1                        &                       1                        &                  $D_j$                  \\ \hline
        $D_{j-1}$                                    &          $D_{j-1}$           &                      1                       &                      1                       &                $D_j$                 &                       1                        &                       1                        &                  $D_j$                  \\ \hline
        $D_{j-1}$                                    &           $D_{t}$            &                    $a_t$                     &                    $a_t$                     &                $D_t$                 &                       1                        &                       1                        &                  $D_j$                  \\ \hline
        $D_{t}$                                      &           $D_{t}$            &                      1                       &                      1                       &                $D_j$                 &                       1                        &                       1                        &                  $D_j$                  \\ \hline
        $D_{t}$                                      &           $D_{q}$            &                     $ a_q $                      &                    $  a_q $                      &                $D_q$                 &                     $1$                      &                     $1$                      &                  $D_j$                  \\ \hline\hline
    \end{longtable}
\end{center}

Since  $S=\bigcup\limits_{i=0}^j D_i$, then   $D_j, D_{j-1}\subseteq S$, $D_t \subseteq S$, for $0\leq t \leq j-2$, so if $x,y \in S$, from the tables we have  that (ii) $x\wedge y, x \vee y \in S$, and  $x \fh y  \in S$. 
Thus from (i) and (ii) it follows that $S$ is a Heyting subalgebra of $FT_n(r)$. 	Let us see that $G \subseteq S$. If $g \in G$, $h(g) \in \{0 = a_0, a_1, \dots, a_{j-1}, 1\}$. 	If $h(g) = 1$, $g \in P_h$ and from (III), $f(g) = 1$, so, $g \in P_{f}$. Then $g \in P_{h} \cap P_{f} = D_{j} \subseteq S$. 	If $h(g) = a_{j-1}$, $g \in Q_{j-1} - Q_{j}$ and from (III), $f(g) = 1$, that is $g \in P_{f}$. So $g \in (Q_{j-1} - Q_{j}) \cap P_{j-1} = D_{j-1} \subseteq S$. 	If $h(g) = a_t$, $0 \leq t \leq j-2$, then $g \in Q_t - Q_{t+1}$ and from (II), $f(g) = a_t$, that is, $g \in P_t - P_{t+1}$. Then $g \in (Q_t - Q_{t+1}) \cap (P_t - P_{t+1}) = D_t \subseteq S$. 	Therefore, $G \subseteq S$. We then have $S = FT_n(r)$.

Then we can write $P_h= Q_{j} = Q_{j} \cap FT_n(r) = Q_{j} \cap \bigcup\limits_{i=0}^{j} D_i = \bigcup\limits_{i=0}^{j}(Q_{j} \cap  D_i) = Q_{j} \cap P_{j-1} =  P_h \cap P_f$. That is, $P_h = P_h \cap P_f$, and therefore $P_h \subseteq P_f$. 	If $P_h = P_f$ then $Ker(\bar{h}) = Ker(\bar{f})$ and then $\bar{h} = \bar{f}$ and $h = f$, which contradicts (IV). Therefore $P_h \subset P_f$.
\end{proof}

\begin{lemma} \label{L56} $|X_M(n,r)|= 2^r$.
\end{lemma}
\begin{proof}	
	If $Y\subseteq G$, let $f:FT_n(r)\to C_n$ defined by $f(g)=1$ if $g\in Y$ and 
	$f(g)=0$ if $g\in G\setminus  Y$. $f$ can be extended to a  homomorphism $\bar{f}$ from $FT_n(r)$ to $C_n$ and $HS_{C_n} (\bar{f}(G))=\bar{f}(G)=\{0,1\}= C_n(2)$. As $\{1\}\in X_M(C_2)$, then $Ker(\bar{f})=P_f\in X_M(n,r)$.  So $|X_M(n,r)|\geq 2^r$.
	
	Let $P \in X_M(n,r)$ and  $h$ the natural  epimorphism  from $FT_n(r)$ onto $FT_n(r)/P= C_n(2)$ and  $f=h|_{G}$,  then $f:G\to C_2$ is such that $\overline{f}=h$ and $P=P_h$, so $|X_M(n,r)|= 2^r$.
\end{proof}

The elements of $X_M(n,r)$ are associated to $r$-tuples formed by zeros and ones.

\begin{example} Using the ordering from Theorem \ref{orderFilters}, we can produce the diagram of $X(4,3)$:
\[
			\beginpicture
			\setcoordinatesystem units <8mm,8mm>
			\setplotarea x from 0 to 12, y from 0 to 5
			
			\put {$ \bullet$} [c] at 0 0 
			\put {$ \bullet$} [c] at 0 5 
			\put {$ \bullet$} [c] at 1 0
            \put {$ \bullet$} [c] at 1 1
            \put {$ \bullet$} [c] at 1 4
            \put {$ \bullet$} [c] at 1 5
			\put {$ \bullet$} [c] at 2 0
            \put {$ \bullet$} [c] at 2 4
            \put {$ \bullet$} [c] at 2 5
			\put {$ \bullet$} [c] at 3 0
            \put {$ \bullet$} [c] at 3 4
            \put {$ \bullet$} [c] at 3 5
			\put {$ \bullet$} [c] at 4 0
            \put {$ \bullet$} [c] at 4 1
            \put {$ \bullet$} [c] at 4 4
			\put {$ \bullet$} [c] at 5 0
            \put {$ \bullet$} [c] at 5 3
            \put {$ \bullet$} [c] at 5 4
            \put {$ \bullet$} [c] at 5 5            
			\put {$ \bullet$} [c] at 6 0
            \put {$ \bullet$} [c] at 6 2
            \put {$ \bullet$} [c] at 6 3
            \put {$ \bullet$} [c] at 6 4
			\put {$ \bullet$} [c] at 7 0
            \put {$ \bullet$} [c] at 7 1
            \put {$ \bullet$} [c] at 7 4
			\put {$ \bullet$} [c] at 8 0
            \put {$ \bullet$} [c] at 8 3
            \put {$ \bullet$} [c] at 8 4
            \put {$ \bullet$} [c] at 8 5  
			\put {$ \bullet$} [c] at 9 0
            \put {$ \bullet$} [c] at 9 1
            \put {$ \bullet$} [c] at 9 3
            \put {$ \bullet$} [c] at 9 4
			\put {$ \bullet$} [c] at 10 0
            \put {$ \bullet$} [c] at 10 1
            \put {$ \bullet$} [c] at 10 4
			\put {$ \bullet$} [c] at 11 0
            \put {$ \bullet$} [c] at 11 1
            \put {$ \bullet$} [c] at 11 3
            \put {$ \bullet$} [c] at 11 4
            \put {$ \bullet$} [c] at 11 5  
            \put {$ \bullet$} [c] at 12 1
            \put {$ \bullet$} [c] at 12 3
            \put {$ \bullet$} [c] at 12 4
            
            \setlinear \plot  1 1 6 2  /
			\setlinear \plot  1 0 1 1  /
            \setlinear \plot  0 0 1 1  /
			\setlinear \plot  4 1 6 2  /
			\setlinear \plot  2 0 1 1  /
			\setlinear \plot  3 0 4 1  /
            \setlinear \plot 4 0 4 1 /
            \setlinear \plot 5 0 4 1 /
            \setlinear \plot 6 0 7 1 /
            \setlinear \plot 7 0 7 1 /
            \setlinear \plot 8 0 7 1 /
            \setlinear \plot 9 0 9 1 /
            \setlinear \plot 10 0 10 1 /
            \setlinear \plot 11 0 11 1 /
            \setlinear \plot 12 1 6 2 /
            \setlinear \plot 4 1 6 2 /
            \setlinear \plot 7 1 6 2 /
            \setlinear \plot 9 1 6 2 /
            \setlinear \plot 10 1 6 2 /
            \setlinear \plot 11 1 6 2 /

			\setlinear \plot  1 4 1 5  /
			\setlinear \plot  2 4 2 5  /
			\setlinear \plot  3 4 3 5  /
            \setlinear \plot 4 4 5 5 /
            \setlinear \plot 5 4 5 5 /
            \setlinear \plot 6 4 5 5 /
            \setlinear \plot 7 4 8 5 /
            \setlinear \plot 8 4 8 5 /
            \setlinear \plot 9 4 8 5 /
            \setlinear \plot 10 4 11 5 /
            \setlinear \plot 11 4 11 5 /
            \setlinear \plot 12 4 11 5 /
            
            \setlinear \plot 5 3 5 4 /
            \setlinear \plot 6 3 6 4 /
            \setlinear \plot 8 3 8 4 /
            \setlinear \plot 9 3 9 4 /
            \setlinear \plot 11 3 11 4 /
            \setlinear \plot 12 3 12 4 /

			\put {\tiny 27} [c] at  .5 1
			\put {\tiny 1} [c] at -.5 5 
			\put {\tiny 29} [c] at 1.4 0
            \put {\tiny 28} [c] at 0.4 0
            \put {\tiny 3} [c] at 1.4 4
            \put {\tiny 2} [c] at 1.4 5
			\put {\tiny 31} [c] at 4.4 1
            \put {\tiny 30} [c] at 2.4 0
            \put {\tiny 5} [c] at 2.4 4
            \put {\tiny 4} [c] at 2.4 5
			\put {\tiny 32} [c] at 3.4 0
            
            \put {\tiny 7} [c] at 3.4 4
            \put {\tiny 6} [c] at 3.4 5
			\put {\tiny 33} [c] at 4.4 0
            \put {\tiny 9} [c] at 4.4 4
			\put {\tiny 34} [c] at 5.4 0
            \put {\tiny 35} [c] at 7.4 1
            \put {\tiny 11} [c] at 5.4 3
            \put {\tiny 10} [c] at 5.4 4
            \put {\tiny 8} [c] at 5.4 5            
			\put {\tiny 36} [c] at 6.4 0
            \put {\tiny 26} [c] at 6.6 2.1
            \put {\tiny 13} [c] at 6.4 3
            \put {\tiny 12} [c] at 6.4 4
			\put {\tiny 37} [c] at 7.4 0
            \put {\tiny 15} [c] at 7.4 4
			\put {\tiny 38} [c] at 8.4 0
            \put {\tiny 39} [c] at 9.4 1
            \put {\tiny 17} [c] at 8.4 3
            \put {\tiny 16} [c] at 8.4 4
            \put {\tiny 14} [c] at 8.4 5  
			\put {\tiny 40} [c] at 9.4 0
            \put {\tiny 19} [c] at 9.4 3
            \put {\tiny 18} [c] at 9.4 4
			\put {\tiny 42} [c] at 10.4 0
            \put {\tiny 21} [c] at 10.4 4
			\put {\tiny 44} [c] at 11.4 0
            
            \put {\tiny 41} [c] at 10.4 .9
            \put {\tiny 43} [c] at 11.4 .9
            \put {\tiny 23} [c] at 11.4 3
            \put {\tiny 22} [c] at 11.4 4
            \put {\tiny 20} [c] at 11.4 5  
			\put {\tiny 45} [c] at 12.4 .9
            \put {\tiny 25} [c] at 12.4 3
            \put {\tiny 24} [c] at 12.4 4
			\endpicture
\]

If we write $a_1$ for $a_{1,4}$ and $a_2$ for $a_{2,4}$, understanding that all the elements are in the chain $C_4$, then corresponding functions are detailed in following table:

\medskip

\centering
\begin{tabular}{||c|c||c|c||c|c||c|c||c|c||} \hline\hline
1  & $(0,0,0)$       & 10 & $(0,a_1,1)$     & 19 & $(a_2,0,a_1)$   & 28 & $(a_1,a_2,a_2)$ & 37 & $(a_2,1,a_1)$     \\
2  & $(0,0,1)$       & 11 & $(0,a_1,a_2)$   & 20 & $(1,1,0)$       & 29 & $(a_1,a_2,1)$   & 38 & $(1,a_2,a_1)$     \\
3  & $(0,0,a_1)$     & 12 & $(0,1,a_1)$     & 21 & $(a_1,a_1,0)$   & 30 & $(a_1,1,a_2)$   & 39 & $(a_1,a_1,1)$     \\
4  & $(0,1,0)$       & 13 & $(0,a_2,a_1)$   & 22 & $(1,a_1,0)$     & 31 & $(1,a_1,1)$     & 40 & $(a_1,a_1,a_2)$   \\
5  & $(0,a_1,0)$     & 14 & $(1,0,1)$       & 23 & $(a_2,a_1,0)$   & 32 & $(a_2,a_1,a_2)$ & 41 & $(a_1,1,a_1)$     \\
6  & $(1,0,0)$       & 15 & $(a_1,0,a_1)$   & 24 & $(a_1,1,0)$     & 33 & $(a_2,a_1,1)$   & 42 & $(a_1,a_2,a_1)$   \\
7  & $(a_1,0,0)$     & 16 & $(a_1,0,1)$     & 25 & $(a_1,a_2,0)$   & 34 & $(1,a_1,a_2)$   & 43 & $(1,a_1,a_1)$     \\
8  & $(0,1,1)$       & 17 & $(a_1,0,a_2)$   & 26 & $(1,1,1)$       & 35 & $(1,1,a_1)$     & 44 & $(a_2,a_1,a_1)$   \\
9  & $(0,a_1,a_1)$   & 18 & $(1,0,a_1)$     & 27 & $(a_1,1,1)$     & 36 & $(a_2,a_2,a_1)$ & 45 & $(a_1,a_1,a_1)$   \\
\hline\hline
\end{tabular}
\label{table:X43}
\end{example}
\medskip

Let $(X,\leq_X)$ and $(Y,\leq_Y)$ be posets with disjoint
underlying sets. Their \emph{cardinal sum} (or disjoint sum),
denoted by $X+Y$, is the poset on $X\cup Y$ whose order restricts
to $\leq_X$ on $X$ and to $\leq_Y$ on $Y$, with no additional
comparabilities between $X$ and $Y$. For the ordinal addition of more than two posets we also use the $\sum$ notation. Their \emph{ordinal sum},
denoted by $X\oplus Y$, is the poset on $X\cup Y$ obtained in the
same way, but with $x<y$ for every $x\in X$ and $y\in Y$.

\begin{remark}
From the proof of Lemma \ref{L56} it follows that $P_{(0,\ldots,0)}$ is both a maximal and minimal element of $X(n,r)$, so it forms by itself a connected component of the poset $X(n,r)$.
\end{remark}

\begin{remark} \label{R615} If $P \in X_M(n,r)$ let $(P]_X= \{ Q\in X_M(n,r): Q\subseteq P\}$. If $P_1,P_2 \in X_M(n,r)$, $P_1 \neq P_2$, then $ (P_1]_{X}\cap (P_2]_{X}= \emptyset$. Indeed, if $Q \in X(n,r)$ is such that $Q \in  (P_1]_{X}\cap (P_2]_{X}$ then $P_1,P_2\in \textbf{S}(Q)$, contradicts Theorem  \ref{T56}. It follows that $X(n,r)=  \sum\limits_{P\in X_M(n,r)}(P]_X$, and by  Lemma \ref{L56},  $X(n,r)$ is the cardinal sum of $2^r$ posets. 
\end{remark}

\begin{remark}
	In what follows we will identify the functions in $F(n,r)$ with their corresponding elements in $X(n,r)$. Our goal is to determine the structure and number of elements of $X(n,r)$ in order to obtain a description of the free algebra $FT_n(r)$.
\end{remark}

\begin{remark}
	We have already defined in Definition \ref{defF_jnr} the sets $F_j(n,r)$. It is clear from the order given in Theorem \ref{orderFilters}, that the elements in $F_2(n,r)$ are the maximal elements in $X(n,r)$. 
\end{remark}

We can discriminate the elements in $F_2(n,r)$ according to the number of 1's appearing in the image of the functions: for $0\leq k \leq r$, let 
\[ F_2(n,r,k)=\set{f:G\to\set{0,1}| |f^{-1}(1)|=k}   \] 

Then, 
\[F_2(n,r)=\bigsqcup_{k=0}^r F_2(n,r,k). \]
Each set $F_2(n,r,k)$ has $\binom{r}{k}$ elements, and $F_2(n,r)$ has $2^r$ elements.

We are interested in determining the poset $(f]_{X(n,r)}$ for each $f\in F_2(n,r)$.

If $f,g\in F_2(n,r,k)$, for $1\le k\le r$,  then $(f]_{X(n,r)}\cong (g]_{X(n,r)}$, where the order isomorphism is a fixed permutations of the components. We can choose a specific function to represent all the ones in the set $F_2(n,r,k)$, namely:
\[ f(n,r,k)=(\overbrace{\underbrace{0,\ldots,0}_{r-k},\underbrace{1\ldots,1}_k}^r)  \]

Now we can reduce the problem noticing that by condition (II) of Theorem \ref{orderFilters} for $t=0$, if a function $f$ has a zero in one coordinate, all the functions in the poset $(f]_{X(n,r)}$ will be zero in that coordinate. Thus, 
\[  (f(n,r,k)]_{X(n,r)}\cong(f(n,k,k)]_{X(n,k)}  \]
In other words, we can drop the coordinates that take value $0$, and determine the posets $(f(n,k)]_{X(n,k)}$ where $f(n,k)$ is the function with $k$ coordinates, all equal to $1$. Notice that there exists for each $r\ge 1$ a function $f(n,0)$ in $X(n,r)$. This function takes the value $0$ for all $g\in G$, and by Theorem \ref{orderFilters}, 
\begin{equation} \label{fn0}
    (f(n,0)]_{X(n,r)}\cong \textbf{1},
\end{equation} where \textbf{1} is the poset with a single element.

For $n=2$ the chains of prime filters cannot have more than a single element, so we also get that for every $k\ge 1$
\begin{equation}\label{f2k}
(f(2,k)]_{X(2,k)}\cong \textbf{1}.
\end{equation}

If $n\ge 3$, we can consider an element $h\in F_3(n,k)$ such that $f(n,k)$ covers $h$. Then we know that there exists at least one $g\in G$ such that $h(g)=a_{1,n}$.  There are $2^k-1$ such functions. In a similar way to what we did before, we group these functions according to the number of occurrences of $a_{1,n}$ in their image. The downsets associated with functions in the same group are pairwise isomorphic and disjoint, so it is enough to compute just one of them. For this reason, we choose one function in each group as a representative:
 \[ h(n,p,s)=(\overbrace{\underbrace{a_{1,n},\ldots,a_{1,n}}_{s},\underbrace{1\ldots,1}_p}^k)  \]

The next theorem proves that  $ (h(n,p,s)]_{X(n,s+p)}\cong(f(n-1,p)]_{X(n-1,p)}$, where $1\le s\le k$ and $p=k-s$. First, we give an example of this fact:

\begin{example}\label{ExDownsetIso}
Let $n=4$, $r=3$, and consider the element $h=h(4,1,2)=(a_{1,4},1,1)\in X(4,3)$.
Since the image of $h$ is $\{a_{1,4},1\}$, the subalgebra it generates in
$\mathbf{C}_4$ is $\{0,a_{1,4},1\}\cong \mathbf{C}_3$, so $h\in F_3(4,3)$.

By Theorem~\ref{orderFilters}, the elements covered by $h$ are obtained by
choosing a nonempty subset of the coordinates equal to $1$ and replacing
them by $a_{2,4}$ (the new value introduced when passing from $j=3$ to
$j+1=4$), while the coordinate equal to $a_{1,4}$ remains unchanged.
Since $h$ has exactly two coordinates equal to $1$ (the second and the
third), there are $2^2-1=3$ such elements, and no further covering is
possible because $j=n=4$ is the maximal level. Hence the decreasing subset
 $(h]_{X(4,3)}$ is $\{(a_{1,4},1,1),(a_{1,4},a_{2,4},1),(a_{1,4},1,a_{2,4}), (a_{1,4},a_{2,4},a_{2,4})\}$.

Now let $n=3$, $r=2$, and consider $f=f(3,2)=(1,1)\in X(3,2)$. The image of $f$
is $\{1\}$, so the subalgebra it generates is $\{0,1\}\cong \mathbf{C}_2$
and $f\in F_2(3,2)$. By the same theorem, the elements covered by $f$ are
obtained by choosing a nonempty subset of its two coordinates (both equal
to $1$) and replacing them by $a_{1,3}$, the value introduced when passing
from $j=2$ to $j+1=3=n$. This gives $2^2-1=3$ elements, and again no
further covering is possible since $j=n=3$. Hence
$(f]_{X(3,2)}$ has elements $(1,1),(a_{1,3},1),(1,a_{1,3})$, and $(a_{1,3},a_{1,3})$.

Both downsets consist of a maximum element covering three pairwise
incomparable elements, so the map $\alpha\colon (h]_{X(4,3)} \to(f]_{X(3,2)}
$ given by
\[
\alpha(a_{1,4},1,1)=(1,1), \quad
\alpha(a_{1,4},a_{2,4},1)=(a_{1,3},1), \]
\[\alpha(a_{1,4},1,a_{2,4})=(1,a_{1,3}), \quad
\alpha(a_{1,4},a_{2,4},a_{2,4})=(a_{1,3},a_{1,3})
\]
is an order isomorphism. The corresponding Hasse diagrams are shown below:

\begin{center}
\setlength{\unitlength}{7mm}
\begin{picture}(24,4)(-2,0)
\put(4,4){\circle*{0.25}}
\put(1,1){\circle*{0.25}}
\put(4,1){\circle*{0.25}}
\put(7,1){\circle*{0.25}}

\put(4,4){\line(-1,-1){2.8}}
\put(4,4){\line(0,-1){3}}
\put(4,4){\line(1,-1){2.8}}

\put(4.5,3.8){$h=(a_{1,4},1,1)$}
\put(-1.9,-0.2){$(a_{1,4},a_{2,4},a_{2,4})$}
\put(2.6,-0.2){$(a_{1,4},a_{2,4},1)$}
\put(6.6,-0.2){$(a_{1,4},1,a_{2,4})$}

\put(16,4){\circle*{0.25}}
\put(13,1){\circle*{0.25}}
\put(16,1){\circle*{0.25}}
\put(19,1){\circle*{0.25}}

\put(16,4){\line(-1,-1){2.8}}
\put(16,4){\line(0,-1){3}}
\put(16,4){\line(1,-1){2.8}}

\put(16.5,3.8){$f=(1,1)$}
\put(12,-0.2){$(a_{1,3},a_{1,3})$}
\put(15,-0.2){$(a_{1,3},1)$}
\put(18.2,-0.2){$(1,a_{1,3})$}

\end{picture}
\end{center}
\end{example}

\begin{theorem} \label{iso} For all $n\ge 3$, $p\ge 1$, and $s\ge 1$ we can establish an order isomorphism:
	\[  (h(n,p,s)]_{X(n,s+p)}\cong(f(n-1,p)]_{X(n-1,p)}  \]
\end{theorem}
\begin{proof} We are going to define the map $\alpha$ from the previous example in a general form. The first consideration is that $\alpha$ drops the first $s$ coordinates of a function, where $h$ takes the value $a_{1,n}$. To formalize this, let  $G_{s+p}=\set{g_1,g_2,\ldots, g_{s+p}}$, and $G_p=\set{g_{s+1},\ldots,g_{s+p}}$. 
For any $u\in$  $(h(n,p,s)]_{X(n,s+p)}$, and $g\in G_p$ let:
\[\alpha(u)(g) = \begin{cases}
a_{i-1,n-1} & \text{if }u(g)=a_{i,n} \text{ for } 2\le i\le n-2 \\
1 & \text{ if } u(g)=1
\end{cases}\]

Observe that if $u\in(h(n,p,s)]_{X(n,s+p)}$, then $\alpha(u)$ is indeed in $(f(n-1,p)]_{X(n-1,p)}$. Let $u\in(h(n,p,s)]$. Assume that $u\in F_j(n,s+p)$, $3\le j\le n$. Then $u(G_{s+p})\supseteq\{a_{1,n},\ldots,a_{j-2,n}\}$ and $u(G_{p})\supseteq\{a_{2,n},\ldots,a_{j-2,n}\}$, so $\alpha(u(G_{p}))\supseteq\{a_{1,n-1},\ldots,a_{j-3,n-1}\}$. 

Thus, $\alpha(u)\in$  $ F_{j-1}(n-1,p)$. Furthermore, from the definition of $\alpha$, $0$ is not in the image of $\alpha(u)$, and that is enough to prove that $\alpha(u)\le f(n-1,p)$.

For $u, v\in (h(n,p,s)]_{X(n,s+p)}$,  $v$ covers $u$ if and only if the following conditions are fulfilled:

\begin{enumerate}[label=\rm{(\roman*)}]
	\item $v\in F_j(n,s+p)$, $u\in F_{j+1}(n,s+p)$, $3\le j\le n-1$.
	\item For all $g\in G_{s+p}$, $v(g)=a_{t,n}$ iff $u(g)=a_{t,n}$, for $1\le t\le j-2$.
	\item For all $g\in G_{s+p}$, $v(g)=1$ iff $u(g)=1$ or $u(g)=a_{j-1,n}$.
	\item There exists $g\in G_p$ such that $v(g)=1$ and $u(g)=a_{j-1,n}$.
\end{enumerate}

On the other hand, $\alpha(v)$ covers $\alpha(u)$ if and only if:

\begin{enumerate}[label=\rm{(\roman*')}]
	\item $\alpha(v)\in F_{j-1}(n-1,p)$, $\alpha(u)\in F_{j}(n-1,p)$, $3\le j\le n-1$.
	\item For all $g\in G_p$, $\alpha(v)(g)=a_{t,n-1}$ iff $\alpha(u)(g)=a_{t,n-1}$, for $1\le t\le j-3$. 
	\item For all $g\in G_p$, $\alpha(v)(g)=1$ iff $\alpha(u)(g)=1$ or $\alpha(u)(g)=a_{j-2,n-1}$.
	\item There exists $g\in G_p$ such that $\alpha(v)(g)=1$ and $\alpha(u)(g)=a_{j-2,n-1}$.
\end{enumerate}

Assume now that conditions (i)-(iv) hold and we prove:
\begin{enumerate}[label=\rm{(\roman*')}]
	\item We know from $v\in F_j(n,s+p)$ that $v(G_{s+p})$ includes the elements  $a_{2,n}$ to $a_{j-2,n}$. Then  $\alpha(v)(G_p)$ includes the elements $a_{1,n-1}$ to $a_{j-3,n-1}$. Therefore, $\alpha(v)\in F_{j-1}(n-1,p)$. Analogously, $\alpha(u)\in F_j(n-1,p)$.
	\item For each $g\in G_p$, 
    and $1\le t\le j-3$ we have that 
    $\alpha(v)(g)=a_{t,n-1}$ iff $v(g)=a_{t+1,n}$. By (ii) it follows that this is equivalent to $u(g)=a_{t+1,n}$, and this in turn is equivalent to $\alpha(u)(g)=a_{t,n-1}$.
    \item $\alpha(v)(g)=1$ exactly when $v(g)=1$, which happens precisely when either $u(g)=1$ or $u(g)=a_{j-1,n}$, equivalently, when either $\alpha(u)(g)=1$ or $\alpha(u)(g)=a_{j-2,n-1}$.
    \item By (iv), there exists $g\in G_p$ such that $v(g)=1$ and $u(g)=a_{j-1,n}$, and applying $\alpha$, we have that$\alpha(v)(g)=1$ and $\alpha(u)(g)=a_{j-2,n-1}$.
\end{enumerate}

Conversely, suppose that $\alpha(v)$ covers $\alpha(u)$, so (i')-(iv') hold. We prove the conditions:
\begin{enumerate}[label=\rm{(\roman*)}]
	\item From (i'), $\alpha(v)\in F_{j-1}(n-1,p)$, so $\alpha(v)(G_p)\supseteq\{a_{1,n-1},\ldots,a_{j-3,n-1}\}$. Then $v(G_p)\supseteq\{a_{2,n},\ldots,a_{j-2,n}\}$, so $v\in F_{j}(n,s+p)$, since $s\ge 1$ and $v(g)=a_{1,n}$ for any $g\in G_{s+p}\setminus G_p$. Similarly with $\alpha(u)$ we deduce that $u(G_p)\supseteq\{a_{2,n},\ldots,a_{j-1,n}\}$, so $u\in F_{j+1}(n,s+p)$.

	\item The following conditions are equivalent for any $g\in G_{s+p}$ and $2\le t\le j-2$: $v(g)=a_{t,n}$; $\alpha(v)(g)=a_{t-1,n-1}$; $\alpha(u)(g)=\alpha(v)(g)=a_{t-1,n-1}$; $u(g)=a_{t,n}$.  
	\item For $g\in G_p$ we also have the equivalence of the conditions: $v(g)=1$; $\alpha(v)(g)=1$; $\alpha(u)(g)=1$ or $\alpha(u)(g)=a_{j-2,n-1}$; $u(g)=1$ or $u(g)=a_{j-1,n}$.
	\item We take the $g\in G_p$ provided by (iv'), and then $g\in G_p\subseteq G_{s+p}$ verifies $v(g)=1$ and $u(g)=a_{j-1,n}$.
\end{enumerate}

$\alpha$ is injective: 	let $u,v \in (h(n,p,s)]_{X(n,s+p)}$. Suppose that for all $g \in G_p$, we have $\alpha(u)(g) = \alpha(v)(g)$. By definition of $\alpha$, $u(g) = v(g)$, if $g \in G_p$. Also, $u(g') = v(g') = a_{1,n}$, for all $g' \in G_{s+p}-G_p$, so $u = v$.
	
Let's also see that $\alpha$ is surjective: let $w \in (f(n-1,p)]_{X(n-1,p)}$, $w(G_p) \supseteq \lbrace a_{1,n-1}, \dots , a_{j,n-1} \rbrace$ for some $j, 3\le n$. Let:
\[u (g) = \begin{cases}
		a_{1,n} & \text{ if } g \in G_{s+p}-G_p \\
		a_{i,n} & \text{ if } g \in G_p \text{ and } w(g) = a_{i-1,n-1}, \text{ for } 2 \leq i \leq j+1 \\
		1 & \text{ if } g \in G_p \text{ and } w(g) = 1 \\
	\end{cases}\]
	
	We can easily check that $u \in (h(n,p,s)]_{X(n,s+p)}$ and $\alpha(u) = w$.\qedhere
\end{proof}

\medskip

In what follows, we denote by $+$ or $\sum$ the cardinal sum of posets, and by $\oplus$ the ordinal sum. For a poset $P$ and a natural number $n$, the product $nP$ means the cardinal sum of $n$ copies of $P$.

From the previous considerations we can now write:
\begin{equation} \label{Xnr_01}
	X(n,r)\cong\sum_{k=0}^{r}{\binom{r}{k}}(f(n,r,k)]_{X(n,r)}\cong \mathbf{1}+\sum_{k=1}^{r}{\binom{r}{k}}(f(n,k)]_{X(n,k)}  
\end{equation}

We get for $n\ge 3$ and $k\ge 2$ the recurrence: 
\begin{equation}
(f(n,k)]_{X(n,k)}\cong
\left(\sum_{s=1}^{k}{\binom{k}{s}}(h(n,k-s,s)]_{X(n-1,k-s)}\right)\oplus \mathbf{1}
\end{equation}
This formula is obtained by partitioning the set of maps in \(F_3(n,k)\) according to the number of free generators that are sent to \(a_{1,n}\). The number of such functions when $s$ generators go to $a_{1,n}$ is $\binom{k}{s}$, and the corresponding downsets are all isomorphic to  the downset $(h(n,k-s,s)]_{X(n-1,k-s)}$. These downsets are all covered by $f(n,k)$, which explains the ordinal addition of a single element poset in the formula.

Now using the isomorphisms from Theorem \ref{iso} and  (\ref{fn0}) we can write:
\begin{equation}
 \label{Xnr_05} 
(f(n,k)]_{X(n,k)}
\cong\left(\sum_{s=1}^{k}{\binom{k}{s}}(f(n-1,k-s)]_{X(n-1,k-s)}\right)\oplus \mathbf{1} \end{equation}

In other words, equations \(\eqref{fn0}\) to \(\eqref{Xnr_05}\) combine two ingredients: first, a counting argument based on the number of functions in \(F_3(n,r)\) having a prescribed pattern of values on the generators; second, the fact that all functions with the same pattern determine isomorphic downsets in \(X(n,r)\). The binomial coefficients record how many functions there are in each class, while the poset isomorphisms show that each class contributes an isomorphic poset  piece to the recurrence.

\begin{example}
	Consider the case with $n=5$ and $r=3$. For the function $f(5,3)=(1,1,1)$ we have:
	\[(f(5,3)]_{X(5,3)}\cong\left({\binom{3}{1}}(f(4,2)]_{X(4,2)}+{\binom{3}{2}}(f(4,1)]_{X(4,1)}+\mathbf{1}\right)\oplus \mathbf{1}\]
	where
	\[(f(4,2)]_{X(4,2)}\cong\left({\binom{2}{1}}(f(3,1)]_{X(3,1)}+\mathbf{1}\right)\oplus \mathbf{1}\]
	\[\cong (2\mathbf{2}+\mathbf{1})\oplus \mathbf{1}\]
	and 
	\[(f(4,1)]_{X(4,1)}\cong\mathbf{1}\oplus \mathbf{1}=\mathbf{2} \quad\quad
	\vcenter{
		\hbox{
			\beginpicture
			\setcoordinatesystem units <6mm,6mm>
			\setplotarea x from 0 to 0, y from 0 to 1
			
			\put {$ \bullet$} [c] at 0 0 
			\put {$ \bullet$} [c] at 0 1 
			
			\setlinear \plot  0 1 0 0  /
			\endpicture
	} }
	\]
	
	Here $\mathbf{2}$ is the chain with two elements.	We also have that $(f(3,1)]_{X(3,1)}\cong\mathbf{1}\oplus \mathbf{1}=\mathbf{2}$
	Therefore:
	\[(f(4,2)]_{X(4,2)}\cong (2\mathbf{2}+\mathbf{1})\oplus \mathbf{1}\quad \quad \vcenter{
		\hbox{
			\beginpicture
			\setcoordinatesystem units <6mm,6mm>
			\setplotarea x from 0 to 2, y from 0 to 2
			
			\put {$ \bullet$} [c] at 0 0 
			\put {$ \bullet$} [c] at 0 1 
			
			\put {$ \bullet$} [c] at 1 0 
			\put {$ \bullet$} [c] at 1 1 
			\put {$ \bullet$} [c] at 1 2
			
			\put {$ \bullet$} [c] at 2 1
			
			\setlinear \plot  1 2 0 1  /
			\setlinear \plot  1 2 1 1  /
			\setlinear \plot  1 2 2 1  /
			\setlinear \plot  0 1 0 0  /
			\setlinear \plot  1 1 1 0  /
			\endpicture
	} }\]
	
	and
	\[(f(5,3)]_{X(5,3)}\cong\left(3 ((2\mathbf{2}+\mathbf{1})\oplus \mathbf{1}) +{3} \mathbf{2}+\mathbf{1}\right)\oplus \mathbf{1}\]
	
	With the corresponding diagram:
	\[
	\beginpicture
	\setcoordinatesystem units <6mm,6mm>
	\setplotarea x from 0 to 12, y from 0 to 4
	
	\put {$ \bullet$} [c] at 0 0 
	\put {$ \bullet$} [c] at 0 1 
	
	\put {$ \bullet$} [c] at 1 0 
	\put {$ \bullet$} [c] at 1 1 
	\put {$ \bullet$} [c] at 1 2
	
	\put {$ \bullet$} [c] at 2 1
	\put {$ \bullet$} [c] at 3 0
	\put {$ \bullet$} [c] at 3 1
	\put {$ \bullet$} [c] at 4 0 
	
	\put {$ \bullet$} [c] at 4 1
	\put {$ \bullet$} [c] at 4 2 
	\put {$ \bullet$} [c] at 5 1 
	\put {$ \bullet$} [c] at 6 0 
	\put {$ \bullet$} [c] at 6 1 
	\put {$ \bullet$} [c] at 6 4 
	\put {$ \bullet$} [c] at 7 0
	\put {$ \bullet$} [c] at 7 1 
	\put {$ \bullet$} [c] at 7 2 
	\put {$ \bullet$} [c] at 8 1 
	
	\put {$ \bullet$} [c] at 9 1 
	\put {$ \bullet$} [c] at 9 2
	\put {$ \bullet$} [c] at 10 1
	\put {$ \bullet$} [c] at 10 2
	\put {$ \bullet$} [c] at 11 1
	\put {$ \bullet$} [c] at 11 2
	\put {$ \bullet$} [c] at 12 2	
	
    \setlinear \plot  6 4 1 2  /
	\setlinear \plot  6 4 4 2  /
	\setlinear \plot  6 4 7 2  /
	\setlinear \plot  6 4 9 2  /
	\setlinear \plot  6 4 10 2  /
	\setlinear \plot  6 4 11 2  /
	\setlinear \plot  6 4 12 2  /
	\setlinear \plot  1 2 0 1  /
	\setlinear \plot  1 2 1 1  /
	\setlinear \plot  1 2 2 1  /
	
	\setlinear \plot  4 2 3 1  /
	\setlinear \plot  4 2 4 1  /
	\setlinear \plot  4 2 5 1  /
	
	\setlinear \plot  7 2 6 1  /
	\setlinear \plot  7 2 7 1  /
	\setlinear \plot  7 2 8 1  /

	\setlinear \plot  9 2 9 1  /
	\setlinear \plot  10 2 10 1  /
	\setlinear \plot  11 2 11 1  /

	\setlinear \plot  0 1 0 0  /
	\setlinear \plot  1 1 1 0  /
	\setlinear \plot  3 1 3 0  /
	\setlinear \plot  4 1 4 0  /
	\setlinear \plot  6 1 6 0  /
	\setlinear \plot  7 1 7 0  /
	
	\endpicture \]
\end{example}

\begin{example}
	Consider instead the case with $n=4$ and $r=3$. For the function $f(4,3)=(1,1,1)$ we have:
	\[(f(4,3)]_{X(4,3)}\cong\left({\binom{3}{1}}(f(3,2)]_{X(3,2)}+{\binom{3}{2}}(f(3,1)]_{X(3,1)}+\mathbf{1}\right)\oplus \mathbf{1}\]
	where
	\[(f(3,2)]_{X(3,2)}\cong\left({\binom{2}{1}}(f(2,1)]_{X(2,1)}+\mathbf{1}\right)\oplus \mathbf{1}\]
	and 
	\[(f(2,1)]_{X(2,1)}\cong\mathbf{1} \]
	Therefore:
	\[(f(3,2)]_{X(3,2)}\cong (2\mathbf{1}+\mathbf{1})\oplus \mathbf{1} \cong 3\mathbf{1}\oplus \mathbf{1} \quad \quad \vcenter{
		\hbox{
			\beginpicture
			\setcoordinatesystem units <6mm,6mm>
			\setplotarea x from 0 to 2, y from 0 to 2
			
			\put {$ \bullet$} [c] at 0 1 
			
			\put {$ \bullet$} [c] at 1 1 
			\put {$ \bullet$} [c] at 1 2
			
			\put {$ \bullet$} [c] at 2 1
			
			\setlinear \plot  1 2 0 1  /
			\setlinear \plot  1 2 1 1  /
			\setlinear \plot  1 2 2 1  /
			\endpicture
	} }\]
	
	and
	\[(f(4,3)]_{X(4,3)}\cong\left(3 (3\mathbf{1}\oplus \mathbf{1}) +{3} \mathbf{2}+\mathbf{1}\right)\oplus \mathbf{1}\]
	
	and now the diagram is:
	\[
	\beginpicture
	\setcoordinatesystem units <6mm,6mm>
	\setplotarea x from 0 to 12, y from 0 to 4
	
	\put {$ \bullet$} [c] at 0 1 
	
	\put {$ \bullet$} [c] at 1 1 
	\put {$ \bullet$} [c] at 1 2
	
	\put {$ \bullet$} [c] at 2 1
	\put {$ \bullet$} [c] at 3 1
	
	\put {$ \bullet$} [c] at 4 1
	\put {$ \bullet$} [c] at 4 2 
	\put {$ \bullet$} [c] at 5 1 
	\put {$ \bullet$} [c] at 6 1 
	\put {$ \bullet$} [c] at 6 4 
	\put {$ \bullet$} [c] at 7 1 
	\put {$ \bullet$} [c] at 7 2 
	\put {$ \bullet$} [c] at 8 1 
	
	\put {$ \bullet$} [c] at 9 1 
	\put {$ \bullet$} [c] at 9 2
	\put {$ \bullet$} [c] at 10 1
	\put {$ \bullet$} [c] at 10 2
	\put {$ \bullet$} [c] at 11 1
	\put {$ \bullet$} [c] at 11 2
	\put {$ \bullet$} [c] at 12 2	
	\setlinear \plot  6 4 1 2  /
	\setlinear \plot  6 4 4 2  /
	\setlinear \plot  6 4 7 2  /
	\setlinear \plot  6 4 9 2  /
	\setlinear \plot  6 4 10 2  /
	\setlinear \plot  6 4 11 2  /
	\setlinear \plot  6 4 12 2  /
	\setlinear \plot  1 2 0 1  /
	\setlinear \plot  1 2 1 1  /
	\setlinear \plot  1 2 2 1  /
	
	\setlinear \plot  4 2 3 1  /
	\setlinear \plot  4 2 4 1  /
	\setlinear \plot  4 2 5 1  /
	
	\setlinear \plot  7 2 6 1  /
	\setlinear \plot  7 2 7 1  /
	\setlinear \plot  7 2 8 1  /
	\setlinear \plot  9 2 9 1  /
	\setlinear \plot  10 2 10 1  /
	\setlinear \plot  11 2 11 1  /
	
	\endpicture \]

\end{example}

We can calculate the number of elements in $X(n,r)$, in two different ways: on one hand, it is the number of surjective functions from a set with $r$ elements to a set of generators of a chain with $n$ or less elements.

If $SF(n,m)$ is the number of functions from a set with $n$ elements onto a set with $m$ elements, it is well known that: 

\begin{equation}\label{surjective} 
	|SF(n,m)| =  \begin{cases}\sum\limits_{i=0}^{m-1} (-1)^i{\binom{m}{i}} (m-i)^n \; \; \; {\rm if} \; \; n\geq m \cr \cr 0  \hspace{4cm}{\rm if}\;\; n < m. \end{cases}
\end{equation} 

If $j=2$, it is clear that $|F_2(n,r)|=2^r$. 

If $2 < j\le n$, the functions in $F_j(n,r)$ of definition \ref{defF_jnr} are the ones that verify \linebreak{$f(G)\cup \{0,1\}= C_n(j)=\{0=a_0,a_1, \cdots,a_{j-2}, 1\}$}, i.e., verify some of the following conditions:
\begin{itemize} 
	\item $f(G)= \{a_1,a_2,\cdots,a_{j-2}\}=X$,
	\item  $f(G)= X \cup \{0\}$,
	\item $f(G)= X \cup \{1\}$,
	\item $f(G)=  C_n(j)$.
\end{itemize}

For $2\leq j \leq n$, taking into consideration the cases above we have that $|F_j(n,r)| =$ \linebreak$ |SF(r,j-2)| + 2 |SF(r,j-1)|+ |SF(r,j)|$. This allows us to determine the number of elements in each connected component of the poset of prime filters of the free algebra, as follows:

\begin{equation}  \label{sumFnr} 
	|F(n,r)|=|X(n,r)|=\sum_{j=2}^{n}|F_j(n,r)|
\end{equation}
\begin{equation}
	|F_j(n,r)|= |SF(r,j-2)| + 2 |SF(r,j-1)|+ |SF(r,j)| \text{ for }2\leq j \leq n
\end{equation}

The other way of calculating the number $|X(n,r)|$ is through the recurrences of equations \(\eqref{fn0}\) to \(\eqref{Xnr_05}\):

\begin{equation} \label{nXnr_01}
	|X(n,r)|=\sum_{k=0}^{r}{\binom{r}{k}}|(f(n,r,k)]_{X(n,r)}|= 1+\sum_{k=1}^{r}{\binom{r}{k}}|(f(n,k)]_{X(n,k)}|  
\end{equation}
for $n\ge 2$, with the initial conditions:
\begin{equation} \label{nXnr_02}
	|(f(2,k)]_{X(2,k)}|= {1}\quad |(f(n,0)]_{X(n,r)}|={1},\end{equation}
and the recurrence: 
\begin{equation} \label{nXnr_03} |(f(n,k)]_{X(n,k)}|=\left(\sum_{s=1}^{k}{\binom{k}{s}}|(f(n-1,k-s)]_{X(n-1,k-s)}|\right)+1 \end{equation}

The equations (\ref{sumFnr}) and (\ref{nXnr_01}) give the same numbers:

\begin{table}[ht]
	\centering
	\caption{Number of elements in $X(n,r)$}
	\label{tab:X_elements}
	\begin{tabular}{||c||ccccc||}
		\hline\hline
		$r \backslash n$ & 2  &  3  & 4  & 5  & 6 \\ \hline\hline
		1         & 2  &  3  & 3  & 3  & 3  \\
		2         & 4  &  9  & 11 & 11 & 11 \\
		3         & 8  & 27  & 45 & 51 &  51  \\
		4         & 16 & 81  & 191   &  275  &299    \\
		5         & 32 & 243 & 813   &  1563  &  2043  \\
		6         & 64 & 729 & 3431   &  8891  &  14771  \\ \hline\hline
	\end{tabular}
\end{table}

Furthermore, knowing the poset $X(n,r)$ we can determine the structure of $FT_n(r)$ and calculate its number of elements:

\begin{remark}
	\label{O16}
	Recall that Priestley duality provides a contravariant equivalence between the category of bounded distributive lattices and lattice homomorphisms and the category of Priestley spaces and continuous order-preserving maps. We denote by $\sf X$ the functor that assigns to each bounded distributive lattice $L$ its Priestley dual space ${\sf X}(L)$ of prime filters of $L$, ordered by inclusion, with the topology having basis sets of the form $\{ x \in X(L) \mid a \in x \}$ for $a \in L$; and by $\sf D$ the functor that assigns to each Priestley space $P$ its bounded distributive lattice ${\sf D}(P)$ of clopen up-sets.
	
	If the Priestley space dual to a bounded distributive lattice has more than one connected component, then the lattice is isomorphic to the direct product of the bounded distributive lattices dual to each connected component of the space. More precisely, for bounded distributive lattices $L_1$ and $L_2$, we have ${\sf X}(L_1 \times L_2) \cong {\sf X}(L_1) + {\sf X}(L_2)$, where $+$ denotes the disjoint union equipped with the coproduct topology and order; conversely, for Priestley spaces $P_1$ and $P_2$, ${\sf D}(P_1 + P_2) \cong {\sf D}(P_1) \times {\sf D}(P_2)$.
\end{remark}

\begin{theorem}
    The lattice structure of the free $T_n$-algebras $FT_n(r)$ are given by the formulas:
    \begin{equation}
\begin{aligned} \label{FTnr_01}
	FT_n(r)=\fD(X(n,r)) & \cong \fD( \mathbf{1}+\sum_{k=1}^{r}{\binom{r}{k}}(f(n,k)]_{X(n,k)})\\
                & \cong\mathbf{2}\times\prod_{k=1}^{r}\fD((f(n,k)]_{X(n,k)})^{\binom{r}{k}}
\end{aligned}
\end{equation}
and 
\begin{equation}
\begin{aligned} \label{FTnr_02}
	\fD((f(n,k)]_{X(n,k)}) & \cong \fD\left(\left(\sum_{s=1}^{k-1}{\binom{k}{s}}(f(n-1,k-s)]_{X(n-1,k-s)}+\mathbf{1}\right)\oplus \mathbf{1}\right)\\
    & \cong \mathbf{1}\oplus (\mathbf{2}\times\prod_{s=1}^{k-1}\fD((f(n-1,k-s)]_{X(n-1,k-s)})^{\binom{k}{s}}),
\end{aligned}
\end{equation}
taking into account also equations \eqref{fn0} and \eqref{f2k}.
\end{theorem}
\begin{proof}
    The result follows from the formulas \eqref{Xnr_01} and \eqref{Xnr_05}, using also the fact that for an element $x$ in a poset, $\fD((x]) \cong \mathbf{1}\oplus\fD((x]\setminus \{x\})$.
\end{proof}

\begin{example} \label{exampleFT43}
    We use the previous theorem and calculate the lattice structure of $FT_4(3)$:
\[FT_4(3)=\fD(X(4,3))\cong\boldsymbol{2}\times\fD((f(4,1)])^3\times\fD((f(4,2)])^3\times\fD((f(4,3)])\]
 We compute the intermediate steps in the recursion:
 \[\fD((f(4,1)])\cong\boldsymbol{1}\oplus\fD((f(3,0)])\cong\boldsymbol{1}\oplus\boldsymbol{2}\cong\boldsymbol{3}\]
\[\fD((f(4,2)])\cong\boldsymbol{1}\oplus(\fD((f(3,1)])^2\times \fD((f(3,0)]))\]
\[\fD((f(3,1)])\cong\boldsymbol{1}\oplus\fD((f(2,0)])\cong\boldsymbol{3}\]
 \[\fD((f(4,3)])\cong\boldsymbol{1}\oplus(\fD((f(3,2)])^3\times\fD((f(3,1)])^3\times\fD((f(3,0)]))\]
\[\fD((f(3,2)])=\boldsymbol{1}\oplus(\fD((f(2,1)])^2\times\fD((f(2,0)])\cong\boldsymbol{1}\oplus\boldsymbol{2}^3\]
So we arrive at 
\[FT_4(3)\cong\boldsymbol{2}\times\boldsymbol{3}^2\times(\boldsymbol{1}\oplus(\boldsymbol{3}^2\times\boldsymbol{2}))^3\times(\boldsymbol{1}\oplus((\boldsymbol{1}\oplus\boldsymbol{2}^3)\times\boldsymbol{3}^2\times\boldsymbol{2}))\]
 \end{example}

\begin{example} If $n\ge 3$, we have that $X(n,1)$ has three elements, and $FT_n(1)$ has six elements:
\noindent $FT_n(1) \cong C_2\times C_3.$

\[
\hbox{
			\beginpicture
			\setcoordinatesystem units <6mm,6mm>
			\setplotarea x from -1 to 7, y from -3 to 1
			
			\put {$ \bullet$} [c] at 0 0 
			\put {$ \bullet$} [c] at 0 1 
			\put {$ \bullet$} [c] at -1 1 

\put {$ \bullet$} [c] at 4 -1 
\put {$ \bullet$} [c] at 5 -2 
\put {$ \bullet$} [c] at 5 0 
\put {$ \bullet$} [c] at 6 -1 
\put {$ \bullet$} [c] at 6 1 
\put {$ \bullet$} [c] at 7 0 

			\setlinear \plot  4 -1 6 1  /
			\setlinear \plot  5 -2 7 0  /
			\setlinear \plot  5 -2 4 -1  /
			\setlinear \plot  6 -1 5 0  /
			\setlinear \plot  7 0 6 1  /
            
			\setlinear \plot  0 1 0 0  /

    \put {\tiny\rm (0)} [c] at  -1.5 1
    \put {\tiny\rm ($a_{1,n}$)} [c] at  .7 0
    \put {\tiny\rm (1)} [c] at  .5 1
            
    \put {$X(n,1)$} [c] at  -.5 -2
    \put {$FT_n(1)$} [c] at  5.5 -3
            
            \endpicture
	}
\]
\end{example}

\begin{example} In the same way, we calculate $X(4,2)$ and $FT_4(2)$:
\[
			\beginpicture
			\setcoordinatesystem units <6mm,6mm>
			\setplotarea x from 0 to 6, y from 3 to 5
			
			\put {$ \bullet$} [c] at 0 5 
            \put {$ \bullet$} [c] at 1 4
            \put {$ \bullet$} [c] at 1 5
            \put {$ \bullet$} [c] at 2 4
            \put {$ \bullet$} [c] at 2 5
            \put {$ \bullet$} [c] at 3 4
            \put {$ \bullet$} [c] at 4 3
            \put {$ \bullet$} [c] at 4 4
            \put {$ \bullet$} [c] at 4 5
            \put {$ \bullet$} [c] at 5 3
            \put {$ \bullet$} [c] at 5 4
            
			\setlinear \plot  1 4 1 5  /
			\setlinear \plot  2 4 2 5  /
            \setlinear \plot 3 4 4 5 /
            \setlinear \plot 4 4 4 5 /
            \setlinear \plot 5 4 4 5 /
            
            \setlinear \plot 4 3 4 4 /
            \setlinear \plot 5 3 5 4 /
			\endpicture
\]

\[FT_4(2)\cong \fD(X(4,2))\cong\boldsymbol{2}\times\fD((f(4,1)])^2\times\fD((f(4,2)])\]
 Using the posets calculated in example \ref{exampleFT43} we obtain:
\[FT_4(2)\cong\boldsymbol{2}\times\boldsymbol{3}^2\times(\boldsymbol{1}\oplus(\boldsymbol{3}^2\times\boldsymbol{2}))\]
$$
\begin{minipage}{10cm}
\beginpicture
\setcoordinatesystem units <.2cm, .2cm>
\setplotarea x from 0 to 20, y from -2 to 22
\setlinear \plot  4 0 4 12 /
\setlinear \plot  10 2 10 14  /
\setlinear \plot  0 4  0 16 /
\setlinear \plot  6 6 6 18 /
\setlinear \plot  4 0 10 2 /
\setlinear \plot  0 4 6 6 /
\setlinear \plot  0 16  6  18 /
\setlinear \plot  4 12  10 14  /
\setlinear \plot  2 2 8 4 /
\setlinear \plot  7 1 3  5 /
\setlinear \plot  7 13 3 17 /
\setlinear \plot  2 14 8 16 /
\setlinear \plot  4 0 0 4 /
\setlinear \plot  10 2 6 6 /
\setlinear \plot  4 12 0 16 /
\setlinear \plot  10 14 6 18 /
\setlinear \plot  2 2 2 14 /
\setlinear \plot  7 1  7 13 /
\setlinear \plot  8 4  8  16 /
\setlinear \plot  3 5 3 17 /
\setlinear \plot   5  3 5 15 /
\put {$   \bullet$} [c] at  5 3
\put {$   \bullet$} [c] at  5 15

\put {$   \bullet$} [c] at  8 4 
\put {$   \bullet$} [c] at  8 16
\put {$   \bullet$} [c] at  0 4
\put {$   \bullet$} [c] at  4 0
\put {$   \bullet$} [c] at  4 12
\put {$   \bullet$} [c] at  6 18
\put {$   \bullet$} [c] at  7 1
\put {$   \bullet$} [c] at  7 13
\put {$   \bullet$} [c] at  6 6
\put {$   \bullet$} [c] at   0 16
\put {$   \bullet$} [c] at  3 5
\put {$   \bullet$} [c] at  10 2 
\put {$   \bullet$} [c] at  10 14  
\put {$   \bullet$} [c] at  6 18  
\put {$   \bullet$} [c] at  2 14  
\put {$   \bullet$} [c] at  3  17
\put {$   \bullet$} [c] at  2  2

\put{$\times$} [c] at 25 10
\endpicture 
\end{minipage}
\begin{minipage}{10cm}
\beginpicture
\setcoordinatesystem units <.2cm, .2cm>
\setplotarea x from 0 to 20, y from -4 to 22

\setlinear \plot  24 0 24 12 /
\setlinear \plot  30 2 30 14  /
\setlinear \plot  20 4 20 16 /
\setlinear \plot  26 6 26 18 /
\setlinear \plot  24 0 30 2 /
\setlinear \plot  20 4 26 6 /
\setlinear \plot  20 16  26  18 /
\setlinear \plot  24 12  30 14  /

\setlinear \plot  22 2 28 4 /
\setlinear \plot  27 1 23  5 /
\setlinear \plot  27 13 23 17 /
\setlinear \plot  22 14 28 16 /

\setlinear \plot  24 0 20 4 /
\setlinear \plot  30 2 26 6 /
\setlinear \plot  24 12 20 16 /
\setlinear \plot  30 14 26 18 /

\setlinear \plot  22 2 22 14 /
\setlinear \plot  27 1 27 13 /
\setlinear \plot  28 4  28  16 /
\setlinear \plot  23 5 23 17 /
\setlinear \plot  25  3 25 15 /

\put {$   \bullet$} [c] at  25 3
\put {$   \bullet$} [c] at  25 15

\put {$   \bullet$} [c] at  23 17

\put {$   \bullet$} [c] at  23 5
\put {$   \bullet$} [c] at  28 4
\put {$   \bullet$} [c] at  28 16
\put {$   \bullet$} [c] at  22 2
\put {$   \bullet$} [c] at  22 14
\put {$   \bullet$} [c] at  24 0
\put {$   \bullet$} [c] at  24 12
\put {$   \bullet$} [c] at  26 18
\put {$   \bullet$} [c] at  27 1
\put {$   \bullet$} [c] at  27 13
\put {$   \bullet$} [c] at  26 6

\put {$   \bullet$} [c] at  30 2 
\put {$   \bullet$} [c] at  30 14  
\put {$   \bullet$} [c] at  20 16  
\put {$   \bullet$} [c] at  20 4  
\put {$   \bullet$} [c] at  23  17
\put {$   \bullet$} [c] at  23  5
\put {$   \bullet$} [c] at  24 -4
\setlinear \plot  24 -4 24 0 /

\endpicture 
\end{minipage}
$$
\end{example}

\begin{corollary} The number of elements of the algebra $FT_n(r)$, for $n\ge 2, r\ge 1$ can be calculated through the formulas:
    \begin{equation}
\label{nFTnr_01}
	|FT_n(r)|={2}\times\prod_{k=1}^{r}|\fD((f(n,k)]_{X(n,k)})|^{\binom{r}{k}}
\end{equation}
and 
\begin{equation}
 \label{nFTnr_02}
	|\fD((f(n,k)]_{X(n,k)})| = {1}+ ({2}\prod_{s=1}^{k-1}|\fD((f(n-1,k-s)]_{X(n-1,k-s)}))|^{\binom{k}{s}}),
\end{equation}
where $|(f(n,0)]_{X(n,r)}|={1}$, and $|(f(2,k)]_{X(2,k)}|={1}$.
\end{corollary}
From those equation we obtain the following table:
\begin{table}[ht]
	\centering
	\caption{Number of elements in $FT_n(r)$}
	\label{tab:X_elements2}
	\begin{longtable}{||c||cccc||}
		\hline\hline
		$n \backslash r$ & 1  &  2     & 3  & 4   \\ \hline\hline
		2         & 4  &  16  & 256    & 65536    \\
		3         & 6  &  162 & 5078214&  781253724434352747138 \\
		4         & 6 &   342 &  14580985662 &  436410120559165875011624386976521095403818486   \\
		5         & 6&342   & 137186159382   & 61 digit number    \\
	 \hline\hline
	\end{longtable}
\end{table}

We know by Lemma \ref{L520} that every  $T_n$-algebra is an L-algebra, but by Remark \ref{O42}, not every L-algebras is a  $T_n$-algebra. Now we are going to describe the relation between the free L-algebras and the free $T_n$-algebras.

\begin{lemma} \label{L2} If $r$ and $n$ are positive integers such that $ r+2 \leq n$, then  $FL(r)$ is a $T_n$-Heyting algebra. 
\end{lemma}
\begin{proof}  By Lemma \ref{L56a}, for any $P \in X(FL(r))$, $FL(r)/P$ is a chain  and by Lemma \ref{L57},  $|FL(r)/P|\leq r+2$.  By hypothesis,  $r+2 \leq n$ so $|\mathbf{S}(P)|\leq n$ for all $P \in X(FL(r))$ and therefore by  Theorem \ref{T416},  $FL(r)$ is $T_n$-Heyting algebra.
\end{proof}

\begin{theorem} \label{isomorphicFree} If $r$ and $n$ are positive integers such that $ r+2 \leq n$, then $FT_n(r)$ and $FL(r)$ are isomorphic Heyting algebras. 
\end{theorem}
\begin{proof}  
	Let $G=\{g_1, \cdots, g_r\}$ and  $G'=\{g'_1, \cdots, g'_r\}$ be sets of free generators of $FT_n(r)$ and $FL(r)$ respectively.  Let $f:G \to FL(r)$ defined by $f(g_i) =g^{'}_i$.  By Lemma \ref{L2} $FL(r)$ is a $T_n$-algebra, so  $f$ extends uniquely to  a Heyting algebra homomorphism   $\overline{f}: FT_n(r) \to FL(r)$.
    
	Since every $T_n$-Heyting algebra is a L-algebra, the map $h:G' \to FT_n(r)$ defined by $h(g'_i)=g_i$ extends uniquely to a Heyting algebra homomorphism   $\overline{h}: FL(r) \to FT_n(r)$.
    
	Then  $\overline{h} \circ \overline{f}$ is an endomorphism of $FT_n(r)$ that acts as the identity on the set $G$ of generators. By the uniqueness of extension from the free generators $G^{'}$ we have $\overline{h} \circ \overline{f}=Id_{FT_n(r)}$. 
	Similarly, the composition  $\overline{f} \circ \overline{h}: FL(r) \to FL(r)$ is the identity on the set $G$ of generators which implies     $\overline{f} \circ \overline{h}=Id_{FL(r)}$. 	Thus $\overline{f}$ and $\overline{h}$ are isomorphisms, proving   $FT_n(r)\cong FL(r)$. 
 \end{proof}

\begin{corollary}  If $r$ and $n$ are positive integers such that $ r+2 \leq n$ then $X(FT_n(r))$ and $X(FL(r))$ are isomorphic posets. 
\end{corollary}

\begin{corollary}
     If $r + 2 \le n \le m$, then $FT_n(r) \cong FT_m(r) \cong FL(r)$. In particular, the
isomorphism type of $FT_n(r)$ does not depend on $n$ for $n \ge r + 2$.
\end{corollary}

\subsubsection*{Use of artificial intelligence tools disclaimer} Two proofs in this paper were obtained with the assistance of large language model tools. The implication $(T_n)\Rightarrow(E_n)$ in Theorem~\ref{T35}, which had resisted a purely algebraic treatment, was found by Claude Sonnet~4.6 (Anthropic); the proof of Theorem~\ref{isomorphicFree} was suggested by Gemini (Google). In both cases the result itself had been conjectured by the authors beforehand, the AI-generated argument was carefully verified, its exposition was edited and rewritten to meet the standards of the paper, and the authors take full responsibility for the correctness and inclusion of both proofs.


\begin{thebibliography}{Mon74b}
	
	\bibitem[AM87]{AbMo}
	Manuel Abad and Luiz Monteiro.
	\newblock On free {$L$}-algebras.
	\newblock In {\em Notes on {M}athematical {L}ogic. {N}o. 34}, volume~35 of {\em Notas L\'{o}gica Matemática}, pages No. 2, 20. Universidad Nacional del Sur, Bah\'{\i}a Blanca, 1987.
	
	\bibitem[Bir35]{Birkhoff1935structure}
	Garrett Birkhoff.
	\newblock On the structure of abstract algebras.
	\newblock {\em Mathematical Proceedings of the Cambridge Philosophical Society}, 31(4):433–454, 1935.
	
	\bibitem[Car86]{Carnielli85meeting}
	Walter Carnielli.
	\newblock Meeting of the association for symbolic logic: Campinas, brazil, 1985.
	\newblock {\em The Journal of Symbolic Logic}, 51(4):1093--1103, 1986.
	
	\bibitem[Hey30]{heyting30}
	A.~Heyting.
	\newblock Die formalen {Regeln} der intuitionistischen {Logik}. {I}, {II}, {III}.
	\newblock {\em Sitzungsber. Preu{{\ss}}. Akad. Wiss., Phys.-Math. Kl.}, 1930:42--56, 57--71, 158--169, 1930.
	
	\bibitem[HK72]{hecht72equational}
	T.~Hecht and Tibor Katri\v{n}\'{a}k.
	\newblock Equational classes of relative {S}tone algebras.
	\newblock {\em Notre Dame J. Formal Logic}, 13:248--254, 1972.
	
	\bibitem[Hor69a]{Hor1}
	Alfred Horn.
	\newblock Free {$L$}-algebras.
	\newblock {\em J. Symbolic Logic}, 34:475--480, 1969.
	
	\bibitem[Hor69b]{Hor}
	Alfred Horn.
	\newblock Logic with truth values in a linearly ordered {H}eyting algebra.
	\newblock {\em J. Symbolic Logic}, 34:395--408, 1969.
	
	\bibitem[Mon54]{AM54}
	Ant\'{o}nio~A. Monteiro.
	\newblock L'arithm\'{e}tique des filtres et les espaces topologiques.
	\newblock In {\em Segundo symposium sobre algunos problemas matem\'{a}ticos que se est\'{a}n estudiando en {L}atino {A}m\'{e}rica, {J}ulio, 1954}, pages 129--162. Centro de Cooperaci\'{o}n Cient\'{\i}fica de la UNESCO para Am\'{e}rica Latina, Montevideo, 1954.
	
	\bibitem[Mon55]{AM55}
	Ant\'onio Monteiro.
	\newblock Axiomes ind\'ependants pour les alg\`ebres de {B}rouwer.
	\newblock {\em Rev. Un. Mat. Argentina}, 17:149--160 (1956), 1955.
	\newblock Also published in \cite{monteiro07works}, pp. 1131-1142.
	
	\bibitem[Mon59]{AM59}
	A.~Monteiro.
	\newblock {\em Curso de {\'A}lgebra de la L{\'o}gica I}.
	\newblock Instituto de Matem{\'a}tica. Universidad Nacional del Sur, Bah{\'i}a Blanca, 1959.
	
	\bibitem[Mon61]{AM61}
	A.~Monteiro.
	\newblock {\em Seminario de L{\'o}gica Matem{\'a}tica}.
	\newblock Instituto de Matem{\'a}tica. Universidad Nacional del Sur, Bah{\'i}a Blanca, 1961.
	
	\bibitem[Mon62]{AM62a}
	Antonio Monteiro.
	\newblock Lin\'{e}arisation de la logique positive de {H}ilbert-{B}ernays.
	\newblock {\em Rev. Un. Mat. Argentina}, 20:308--309, 1962.
	
	\bibitem[Mon70]{LM70}
	Luiz Monteiro.
	\newblock Les alg\`ebres de {H}eyting et de {{\L}}ukasiewicz trivalentes.
	\newblock {\em Notre Dame J. Formal Logic}, 11:453--466, 1970.
	
	\bibitem[Mon72]{LM72}
	Luiz Monteiro.
	\newblock Alg\`ebre du calcul propositionnel trivalent de {H}eyting.
	\newblock {\em Fund. Math.}, 74(2):99--109, 1972.
	
	\bibitem[Mon74a]{AM74}
	A.~Monteiro.
	\newblock L'arithm{\'e}tique des filtres et les espaces topologiques, 1974.
	\newblock Memoria presentada en 1950 a la Sociedad Matem{\'a}tica de Francia con motivo de la jubilaci{\'o}n de M. Fr{\'e}chet. Publicada en 1974.
	
	\bibitem[Mon74b]{AM74a}
	A.~Monteiro.
	\newblock {\em L'arithm{\'e}tique des filtres et les espaces topologiques I--II}.
	\newblock Number 29--30 in Notas de L{\'o}gica Matem{\'a}tica. Instituto de Matem{\'a}tica, Universidad Nacional del Sur, Bah{\'i}a Blanca, Argentina, 1974.
	
	\bibitem[Mon77]{LM77}
	Luiz Monteiro.
	\newblock Alg\`ebres de {H}ilbert {$n$}-valentes.
	\newblock {\em Portugal. Math.}, 36(3-4):159--173 (1980), 1977.
	
	\bibitem[Mon80]{AM80}
	Ant\'onio Monteiro.
	\newblock Sur les alg\`ebres de {H}eyting sym\'etriques.
	\newblock {\em Portugaliae Mathematica}, 1-4(39):1---237, 1980.
	\newblock Also published in \cite{monteiro07works}, pp. 1769-2037.
	
	\bibitem[Mon95]{AM95b}
	Ant{\'o}nio Monteiro.
	\newblock {\'A}l\-ge\-bras de {H}ey\-ting.
	\newblock Informe T\'ecnico Interno~51, Instituto de Matem\'atica, Universidad Nacional del Sur, Bah\'\i a Blanca, 1995.
	\newblock Also published in \cite{monteiro07works}, pp. 2061-2096.
	
	\bibitem[Mon07]{monteiro07works}
	Ant\'onio~A. Monteiro.
	\newblock {\em The works of Ant\'onio A. Monteiro}.
	\newblock Funda\c{c}\~ao Calouste Gulbenkian, Lisbon and The Humboldt Press, London, 2007.
	
	\bibitem[Sci26]{Scirica26Heyting}
	Carlos~Ernesto Scirica.
	\newblock {\em Álgebras de Heyting de interiores}.
	\newblock PhD thesis, Universidad de Buenos Aires, Ciudad de Buenos Aires, Argentina, Febrero 2026.
	
	\bibitem[Tho62]{TI}
	Ivo Thomas.
	\newblock Finite limitations on {D}ummet's {LC}.
	\newblock {\em Notre Dame J. Formal Logic}, 3:170--174, 1962.
	
	\bibitem[Var66]{Var}
	Jules Varlet.
	\newblock On the characterization of {S}tone lattices.
	\newblock {\em Acta Sci. Math. (Szeged)}, 27:81--84, 1966.
	
\end{thebibliography}
\end{document}